%% file: 2005-41.tex
\documentclass{gtart_h}  

\input gtoutput

\lognumber{561}
\received{5 February 2005}
\volumenumber{9}\papernumber{41}\volumeyear{2005}
\pagenumbers{1835}{1880}   
\revised{11 September 2005}
\published{25 September 2005}
\accepted{7 August 2005}
\proposed{Benson Farb}
\seconded{Martin Bridson, Joan Birman}

\usepackage{amsmath,amssymb,graphicx,color,dcpic}

\def\fref#1{\hyperlink{#1anchor}{\ref*{#1}}}
\def\figref#1{\hyperlink{#1anchor}{Figure~\ref*{#1}}}
\def\anchor#1{\noindent\hypertarget{#1anchor}{\smash{$\phantom{99}$}}\newline}

\newtheorem{theorem}{Theorem}[section]
\newtheorem{lemma}[theorem]{Lemma}
\newtheorem{corollary}[theorem]{Corollary}
\newtheorem{proposition}[theorem]{Proposition}
\theoremstyle{definition}
\newtheorem{notation}[theorem]{Notation}
\newtheorem{definition}[theorem]{Definition}
\newtheorem{remark}[theorem]{Remark}
\newtheorem{example}[theorem]{Example}
\newtheorem{question}{Question}

\DeclareMathOperator{\minlex}{\mathsf{minlex}}
\DeclareMathOperator{\complexity}{\mathsf{c}}
\DeclareMathOperator{\lex}{\mathsf{lex}}
\DeclareMathOperator{\tight}{\mathsf{tight}}
\DeclareMathOperator{\core}{\mathsf{core}}
\DeclareMathOperator{\collapse}{\mathsf{collapse}}
\DeclareMathOperator{\op}{op}
\DeclareMathOperator{\id}{Id}
\DeclareMathOperator{\Aut}{\mathsf{Aut}}

\def\abs(#1,#2){|#2|_{{#1}}}
\def\treeproj(#1){{q_{#1}}}
\def\result(#1){{#1'}}
\def\induce(#1){{\hat{#1}}}

\newcommand{\G}{{\mathcal G}}
\def\bsgraph{\Gamma}
\def\bse{{e}} 
\def\bseb{{f}} 
\def\bsv{{v}} 
\def\bsvb{{u}} 
\def\output{{out}}
\def\splitting{cleaving}
\def\Splitting{Cleaving}
\def\splited{cleaved}

\def\stallings{{\mathsf{S}}}
\def\gersten{{\mathsf{G}}}
\def\graph{\Sigma}
\newcommand{\edge}{e}
\def\e{c} 

\newcommand{\graphseq}{\vec\graph}
\newcommand{\map}{g}
\newcommand{\edgepath}{\sigma}
\newcommand{\edgeset}{\mathcal E}
\newcommand{\morphism}{\map}
\newcommand{\es}{\edge} 
\newcommand{\et}{f} 
\newcommand{\proper}{strict}
\def\proj(#1){{\lambda_{#1}}} 

\def\f(#1){{F_{#1}}} 
\newcommand{\trivialgroup}{{\mathbf 1}}
\newcommand{\group}{G}
\newcommand{\groupo}{G_0}
\newcommand{\groupi}{G'}
\newcommand{\groupii}{G''}
\newcommand{\subgroup}{H}
\newcommand{\conjugacy}{\mathcal H}
\newcommand{\conjugacyo}{\mathcal H_0}
\newcommand{\subgroupo}{H_0}
\newcommand{\subgroupoi}{H_0'}
\newcommand{\subgroupoii}{H_0''}

\newcommand{\groupseq}{\vec H}
\newcommand{\conjugacyseq}{{\vec{\mathcal H}}}
\newcommand{\conjugate}{\sim}
\newcommand{\isomorphic}{\cong} 

\newcommand{\B}{{\mathcal B}}
\newcommand{\basis}{\B}
\newcommand{\word}{w}
\newcommand{\lab}{b}
\newcommand{\llab}{c}

\newcommand{\auto}{\alpha}

\newcommand{\Z}{\mathbb Z}

\def\t{{t}}
\def\s{{s}}

\newcommand{\D}{{\mathcal D}}

\newcommand{\lexspace}{\mathcal L}

\begin{document}

\title[The Grushko decomposition of a
finite graph of free groups]{The Grushko decomposition of a
finite graph\\of finite rank free groups: an algorithm}

\authors{Guo-An Diao\\Mark Feighn}

\address{School of Arts and Sciences, Holy Family 
University\\Philadelphia, PA 19114, USA\\{\rm and}\\Department of Mathematics 
and Computer Science\\Rutgers University, Newark, NJ 07102, USA}
\asciiaddress{School of Arts and Sciences, Holy Family 
University\\Philadelphia, PA 19114, USA\\and\\Department of Mathematics 
and Computer Science\\Rutgers University, Newark, NJ 07102, USA}

\asciiemail{gdiao@holyfamily.edu, feighn@andromeda.rutgers.edu}
\gtemail{\mailto{gdiao@holyfamily.edu}{\rm\qua 
and\qua}\mailto{feighn@andromeda.rutgers.edu}}

\begin{abstract} 
A finitely generated group admits a decomposition, called its {\it
Grushko decomposition}, into a free product of freely indecomposable
groups. There is an algorithm to construct the Grushko decomposition
of a finite graph of finite rank free groups. In particular, it is
possible to decide if such a group is free.
\end{abstract}
\asciiabstract{%
A finitely generated group admits a decomposition, called its Grushko
decomposition, into a free product of freely indecomposable
groups. There is an algorithm to construct the Grushko decomposition
of a finite graph of finite rank free groups. In particular, it is
possible to decide if such a group is free.}




\maketitlepage

\section{Introduction}
\begin{theorem}[Grushko \cite{ig:grushko}]\label{t:grushko}
A finitely generated group $\group$ is a free product of a finite rank
free subgroup and finitely many freely indecomposable non-free subgroups.
\end{theorem}

Up to reordering and conjugation, the non-free factors appearing in
this {\it the Grushko decomposition of $\group$} are unique. The rank
of the free factor is also an invariant of $\group$.  The main result
of this paper is:

\begin{theorem}\label{t:main}
There is an algorithm which produces the Grushko decomposition of a
finite graph\footnote{in the sense of Bass--Serre \cite{se:trees}} of
finite rank free groups\footnote{ie, vertex and edge groups are
  free of finite rank, see Section~\ref{s:our case}}.
\end{theorem}

A relative version is given in Section~\ref{ss:relative}.

The class of finite graphs of finite rank free groups is fascinating
and has received much attention. For example, mapping tori of free
group automorphisms are in this class, see \cite{bf:combination,
fh:coherence, gmsw:hopfian, pb:splittings}. Also, limit groups which
appear in the recent work on the Tarski problem, see \cite{km:limit1,
zs:tarski1}, have a hierarchy in which those limit groups appearing on
the first level are finite graphs of finite rank
free groups.

The algorithm is given in
Section~\ref{s:algorithm}. Theorem~\ref{t:main 2} is a refined version
of Theorem~\ref{t:main} and is proved in Section~\ref{s:final}.  Work
of Shenitzer or Swarup when combined with Whitehead's algorithm for
deciding if a given element of a free group is primitive\footnote{an
element of some basis} gives the case of Theorem~\ref{t:main} where
edge groups are cyclic.

\begin{theorem}[Shenitzer--Swarup \cite{sh:amalg},\cite{swarup:hnn},
\cite{jhcw:equivalent},\cite{jhcw:certain}, see also
\cite{js:grushko}] \label{t:intro}
There is an algorithm to decide whether
a finite graph of finite rank free groups with cyclic edge groups is
free.
\end{theorem}

There are other notable situations where the Grushko decomposition may
be found algorithmically. For example, given a presentation with one
relation for a group $G$, the Grushko decomposition of $G$ can be
constructed (eg, \cite{ls:book}). Also, from a triangulation of a
closed orientable 3--manifold $M$, the connected sum decomposition of
$M$ can be found (see \cite{jlr:algorithms}). The referee points out
two other papers \cite{vg:connectedness} and \cite{km:algorithms} that
consider respectively hyperbolic groups and limit groups.

Here is a sketch of the proof of Theorem~\ref{t:main}. There are three
steps. Suppose that $S$ is a cocompact $\group$--tree with finitely
generated edge stabilizers. Suppose further that $\group$ is freely
decomposable and that $T$ is a $\group$--tree with one orbit of edges
and with trivial edge stabilizers. Give the product $S\times T$ the
diagonal $\group$--action. $S\times T$ is a union of squares
(edge$\times$edge). There is a cocompact $\group$--subcomplex $X_S(T)$
that is a $\group$--deformation retract of $S\times T$ (see
Section~\ref{s:examples}). As $X_S(T)$ is contained in $S\times T$,
there is a natural map $X_S(T)\to T$. The preimage in $X_S(T)$ of the
midpoint of an edge of $T$ is a compact forest. A valence one vertex
of this forest corresponds to a square in $X_S(T)$ that may be
equivariantly collapsed. We may iteratively collapse until each
component of this tree is a point. These points exhibit free
decompositions of $\group$ that are compatible with the original
splitting given by $S$. 

An argument similar to the one in this first step was used by
Bestvina--Feighn in \cite{bf:outerlimits} to among other things
reprove Theorem~\ref{t:intro}.  The use of products of trees is
inspired by the Fujiwara--Papasoglu \cite{fp:jsj} approach to the
Rips--Sela $JSJ$--theorem \cite{rs:jsj}.

The second step is to translate these collapses of $X_S(T)$ into
corresponding {\it simplifications} of the original graph of
groups. This is straightforward and is done in
Section~\ref{s:algebra}. These first two steps do not use the
hypothesis that edge and vertex groups are free.

In the third step (Section~\ref{s:algorithms}), we show how Gersten
representatives \cite{sg:whitehead} of conjugacy classes of subgroups
of free groups can be used to detect simplifications. This is probably
the heart of the paper.

A special case of Theorem~\ref{t:main} solves Problem~F24b on the
problem list at \url{http://www.grouptheory.org}.

For the convenience of readers interested primarily in using the
algorithm, it is described in the next section
(Section~\ref{s:algorithm}). Definitions are given, but proofs are,
for the most part, deferred until later in the paper.  The first two
steps of the proof are more general and therefore somewhat cleaner,
see Sections~\ref{s:square complexes}--\ref{s:algebra} which can be
read independently of Section~\ref{s:algorithm}.

The first named author's thesis \cite{gad:thesis} included an
algorithm to decide if a finite graph of finite rank free groups is
free. The second named author warmly thanks Mladen Bestvina for
helpful conversations and gratefully acknowledges the support of the
National Science Foundation.

\section{The algorithm}\label{s:algorithm}
Mainly to establish notation, we first recall the definition of a graph of groups.
\subsection{Graphs of groups}\label{s:graphs of groups}
A reference for this section is 
\cite{se:trees}. 
A {\it graph} is a 1--dimensional $CW$--complex and is determined by the
following combinatorial data: a 4--tuple $(V,\ \hat E,\ \op,\
\partial_0)$ where
\begin{itemize}
\item
$V$ and $\hat E$ are sets;
\item
$\op\co\hat E\to \hat E$ satisfies
\begin{enumerate}
\item
$\op\circ \op=\id$, and
\item
$\op(\bse)\not=\bse$, for all $\bse\in \hat E$; and
\end{enumerate}
\item
$\partial_0\co \hat E\to V$.
\end{itemize}
For $\bse\in \hat E$, we also write $\bse^{-1}$ for $\op(\bse)$ and
set $\partial_1\bse=\partial_0\bse^{-1}$. For $\bsv\in V$, $\hat
E(\bsv)=\{\bse\in \hat E : \partial_0\bse=\bsv\}$. The {\it valence of
$\bsv$} is the cardinality $|\hat E(\bsv)|$ of $\hat E(\bsv)$.  Such a
4--tuple is a {\it combinatorial graph}.  The graph
$$\bsgraph=\bsgraph(V,\ \hat E,\ \op,\ \partial_0)$$ so determined has
vertex set identified with $V$. The set $\hat E$ corresponds to the
set of oriented edges of $\bsgraph$; the set $E$ of edges of
$\bsgraph$ is identified with $\{\{\bse,\bse^{-1}\}:\bse\in \hat
E\}$. The map $\op$ reverses edge orientations, and $\partial_0$
determines the characteristic maps of $\bsgraph$.  Up to isomorphism,
a graph uniquely determines a combinatorial graph and {\it vice
versa}. In particular, properties of one give properties of the
other. The interior of an edge $\bse$ is denoted $\mathring{\bse}$.

A {\it graph of groups} is a 4--tuple $$\G=(\bsgraph(V,\ \hat E,\ \op,\
\partial_0),\ \{G_{\bsv}:\bsv\in V\},\ \{G_{\bse}:\bse\in \hat E\},\
\{\varphi_{\bse}\co\bse\in \hat E\})$$ where
\begin{itemize}
\item
$\bsgraph(V,\ \hat E,\ \op,\ \partial_0)$ is a connected combinatorial
graph $\bsgraph(\G)$;
\item
for $\bse\in \hat E$ and $\bsv\in V$, $G_{\bse}=G_{\bse^{-1}}$ and
$G_{\bsv}$ are groups; and
\item
for $\bse\in \hat E$, $\varphi_{\bse}\co G_{\bse}\to G_{\partial_0\bse}$
is a monomorphism.
\end{itemize}
The groups $G_\bse$ and $G_\bsv$ are respectively {\it edge} and {\it
  vertex groups}. The $\varphi_{\bse}$'s are {\it bonding maps}.
We say that $\G$ is {\it reduced} if
\begin{itemize}
\item
$\varphi_{\bse}\co G_{\bse}\to G_{\bsv}$ is not an isomorphism for any
valence one vertex $\bsv$; and
\item
if $\bsv$ has valence two and if $\varphi_{\bse}\co G_{\bse}\to G_{\bsv}$ is an
isomorphism, then $\hat E(\bsv)=\{\bse,\bse^{-1}\}$ (in which case
$\bsgraph(\G)$ is a {\it loop}).
\end{itemize}

Associated to a graph of groups $\G$ is an isomorphism type of group
$\pi_1(\G)$, see \cite{se:trees}.  
If $\group\isomorphic\pi_1(\G)$ then we say
that $\G$ is a {\it graph of groups decomposition for $\group$}.
Let $\G$
and $\G'$ be graphs of groups with the same underlying graphs, edge
groups, and vertex groups. We say that $\G$ and $\G'$ are {\it
conjugate}, written $\G\conjugate\G'$, if there is a sequence
$\vec{h}=\{h_{\bse}\in\group_{\partial_0\bse}:\bse\in\hat E\}$ such that $\varphi_{\bse}'=i_{h_{\bse}}
\circ\varphi_{\bse}$ where $i_{h_{\bse}}$ denotes the inner
automorphism induced by $h_{\bse}$, ie, $i_{h_\bse}(g)=h_\bse g h_\bse^{-1}$. If $\G$ and $\G'$ are conjugate, then
$\pi_1(\G')$ and $\pi_1(\G)$ are isomorphic.

A simplicial action of a group $G$ on a tree $T$ determines a graph of
groups with underlying graph $T/G$ and with vertex and edge groups given by
vertex and edge stabilizers of in $T$. Conversely, a graph of groups
$\G$ determines up to simplicial isomorphism a
$\group\cong\pi_1(\G)$--tree $T(\G)$. See \cite{se:trees}. $\G$ is a
{\it trivial} graph of groups decomposition if $\group\cong\pi_1(\G)$
fixes a point of $T(\G)$. If $\G$ is a non-trivial graph of groups
decomposition for $\group$ then we also say that $\G$ is a {\it
  splitting} for $\group$. If the edge groups of $\G$ are contained in
some class of groups then we say that $\group$ {\it splits over this
  class.} For example, if $\G$ is a non-trivial graph of groups
decomposition for $\group$ and all edge groups of $\G$ are trivial then we
say that $\group$ splits over $\trivialgroup$ where $\trivialgroup$
denotes the trivial group.

If $\G$ is a graph of groups with $\group\cong\pi_1(\G)$ then the
edges $\bse$ of $\bsgraph(\G)$ with $G_\bse=\trivialgroup$ determine a
free product decomposition $F_m*G_1*\cdots *G_n$ where $m$ is the rank
of the graph obtained from $\bsgraph(\G)$ by collapsing all edges $\bseb$
with $G_\bseb\not=\trivialgroup$ and the $G_i$'s are the fundamental
groups of graphs of groups given by the components of
$\bsgraph(\G)\setminus(\cup_\bse\{\mathring{\bse}\mid G_\bse=\trivialgroup\}).$
This decomposition is called the {\it decomposition of $\group$ determined by
the edges of $\G$ with trivial stabilizer}.

Now we describe some operations on a graph of groups $\G$. These will
be the simplifying moves of the algorithm. The moves will transform
$$\G=(\bsgraph,\{G_\bsv\},\{G_\bse\},\{\varphi_\bse\})$$ into
$\G'=(\bsgraph',\{G_{\bsv'}\},\{G_{\bse'}\},\{\varphi_{\bse'}\}).$ Much
of the data will be the same for $\G$ and $\G'$, so in describing the
moves we will usually only record the differences.
\subsection{Reducing}
If a bonding map at a valence one or two vertex is an isomorphism then
there is an obvious simplification that we now describe.

Suppose $\bsv\in V$ has valence one, ie, $\hat
E(\bsv)=\{\bse\}$. Suppose further that $\varphi_{\bse}$ is an
isomorphism. Then, define $\G'$ by setting $V'=V\setminus\{\bsv\}$ and
$\hat E'=\hat E\setminus\{\bse,\bse^{-1}\}$.

Next suppose that $\bsv\in V$ has valence two and that $\bsgraph$ is
not a loop, ie, $\hat E(\bsv)=\{\bse,
\bseb\not=\bse^{-1}\}$. Suppose further that $\varphi_{\bse}$ is an
isomorphism. Then, define $\G'$ by setting $V'=V\setminus\{\bsv\}$,
setting $\hat E'=\hat E\setminus\{\bse,\bse^{-1}\}$, redefining $\partial_0\bseb$ to be
$\partial_1\bse$, and redefining $\varphi_{\bseb}$ to be
$\varphi_{\bse^{-1}}\circ\varphi^{-1}_{\bse}\circ\varphi_{\bseb}$.

\begin{figure}[ht!]\anchor{f:reduce}
\cl{\scalebox{0.75}{\input{reduce.pstex_t}}}
\caption{Reducing}\label{f:reduce}
\end{figure}
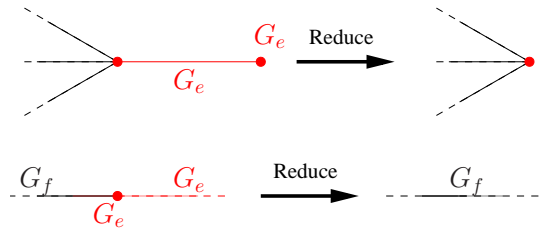

If $\bsgraph(\G)$ is finite then, since each of these operations
decreases the number of vertices, after finitely many operations we
obtain a reduced graph of groups that has been obtained from $\G$ by
{\it reducing}. See \figref{f:reduce}.

In the remaining moves, some vertex group admits a graph of groups
decomposition that is compatible with the bonding maps.

\subsection{Blowing up}\label{s:blowing up}
There are two types of blowing up. Suppose first that, for some $\bsv\in V$,
$\group_{\bsv}=\group_{\bsv}'*\langle t\rangle$ where $t$ has infinite
order, and $\varphi_{\bse}(\group_{\bse})\subset \group_{\bsv}'$ for
$\bse\in \hat E(\bsv)$. Then, define $\G'$ as follows:
\begin{itemize}
\item
the vertex sets are the same, ie, $V'=V$;
\item
add a new oriented loop so that $\hat E'=\hat
  E\cup\{\bse_t,\bse_t^{-1}\}$, with
  $\partial_0\bse_t=\partial_0\bse_t^{-1}=\bsv$,
  $\group_{\bse_t}=\trivialgroup$, and redefine
  $\group_{\bsv}$ to be $\group_{\bsv}'$; and
\item
bonding maps are given by restricting the codomains of the bonding
maps of $\G'$ if necessary.
\end{itemize}
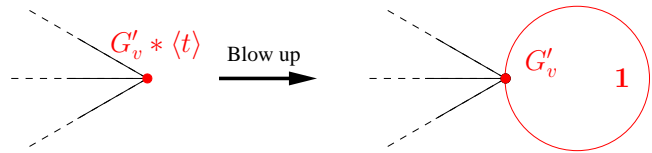
\begin{figure}[ht!]\anchor{f:blowup 1}
\cl{\scalebox{0.75}{\input{blowup2.pstex_t}}}
\caption{The first type of blowup}\label{f:blowup 1}
\end{figure}

Secondly, suppose that, for some $\bsv\in V$,
$\group_{\bsv}=\group_{\bsv}'*\group_{\bsv}''$ and, for each
$\bse\in\hat E(\bsv)$, $\varphi_\bse(\group_\bse)$ is either contained
in $\group_\bsv'$ or $\group_\bsv''$. Then, define $\G'$ as follows:
\begin{itemize}
\item
replace $\bsv$ by two vertices $\bsv'$ and
$\bsv''$, ie, $V'=V\cup\{\bsv',\bsv''\}\setminus\{\bsv\}$;
\item
add a new oriented edge so that $\hat E'=\hat
E\cup\{\bse_t,\bse_t^{-1}\}$ with $\partial_0\bse_t=\bsv'$,
$\partial_1\bse_t=\bsv''$, and $G_{\bse_t}={\trivialgroup}$;
\item
if $\bse\in\hat E(\bsv)$ then
$\partial_0\bse$ is $\bsv'$ or $\bsv''$ depending on whether
$\group_\bse\subset\group_{\bsv'}$ or
$\group_\bse\subset\group_{\bsv''}$; and
\item
bonding maps are given by restricting the codomains of the bonding
maps of $\G'$ if necessary.
\end{itemize}
\begin{figure}[ht!]\anchor{f:blowup2}
\cl{\scalebox{0.75}{\input{blowup3.pstex_t}}}
\caption{The second type of blowup}\label{f:blowup2}
\end{figure}
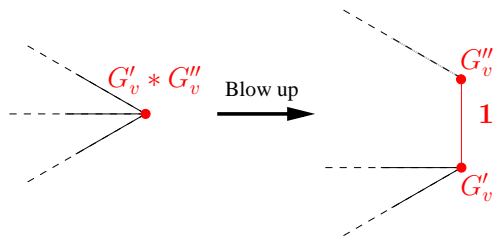
In each of these cases, we say that the new graph of groups is
obtained from $\G$ by {\it blowing up}. See Figures~\fref{f:blowup 1}
and \fref{f:blowup2}.

\subsection{Unpulling}\label{s:unpulling}
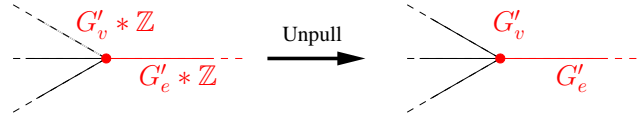
\begin{figure}[ht!]\anchor{f:unpulling}
\cl{\scalebox{0.75}{\input{unpull.pstex_t}}}
\caption{Unpulling}\label{f:unpulling}
\end{figure}

Suppose that, for some $\bsv\in V$ and $\bse\in \hat E(\bsv)$, we have
$\group_{\bsv}=\group_{\bsv}'*\Z$, $\group_{\bse}=\group_{\bse}'*\Z$,
$\varphi_{\bse}(\group_{\bse}')\subset\group_{\bsv}'$,
$\varphi_{\bse}(\Z)=\Z$, and $\varphi_{\bseb}(\group_{\bseb})\subset
\group_{\bsv}'$ for $\bseb\in \hat E(\bsv)\setminus\{\bse\}$. Then,
define $\G'$ as follows:
\begin{itemize}
\item
$V'=V$;
\item
$\hat E'=\hat E$;
\item
redefine $\group_{\bsv}$ to be $\group_{\bsv}'$,
$\group_{\bse}$ to be $\group_{\bse}'$; and
\item
bonding maps are given by restricting codomains of bonding maps of
 $\G'$ if necessary.
\end{itemize}
We say that the new graph of groups is obtained from $\G$ by {\it
unpulling}.\footnote{We use the term {\it unpulling} because the
inverse operation {\it pulls an element of $G_v'$ across the edge
$\bse$}.} See \figref{f:unpulling}.

\subsection{Unkilling}\label{s:unkilling}
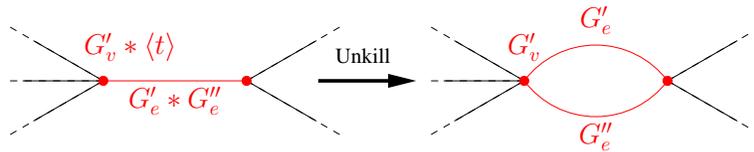
\begin{figure}[ht!]\anchor{f:unkilling}
\cl{\scalebox{0.75}{\input{unkill.pstex_t}}}
\caption{Unkilling}\label{f:unkilling}
\end{figure}

Suppose that, for some $\bsv\in V$ and $\bse\in \hat E(\bsv)$,
$\group_{\bsv}=\group_{\bsv}'*\langle t\rangle$ where $t$ has infinite
order,
$\group_{\bse}=\group_{\bse}'*\group_{\bse}''$,
$\varphi_{\bse}(\group_{\bse}')\subset \group_{\bsv}'$,
$\varphi_{\bse}(\group_{\bse}'')\subset t\group_{\bsv}'t^{-1}$, and
$\varphi_{\bseb}(\group_{\bseb})\subset \group_{\bsv}'$ for $\bseb\in
\hat E(\bsv)\setminus\{\bse\}$. Then, define $\G'$ as follows:
\begin{itemize}
\item
$V'=V$;
\item
the oriented edge $\{\bse,\bse^{-1}\}$ is replaced with two oriented
edges having the same endpoints as $\bse$: $$\hat E'=\hat
E\cup\{\bse',\bse'',(\bse')^{-1},(\bse'')^{-1}\}\setminus\{\bse,\bse^{-1}\}$$
with $\partial_0\bse'=\partial_0\bse''=\partial_0\bse$ and
$\partial_1\bse'=\partial_1\bse''=\partial_1\bse$;
\item
$\group_{\bse'}=\group_{\bse}'$,
$\group_{\bse''}=\group_{\bse}''$, $\group_{\bsv}$ is redefined to be
$\group_{\bsv}'$; and
\item
$\varphi_{\bse''}=i_{t^{-1}}\circ\varphi_{\bse}|_{\group_{\bse}''}$
and other bonding maps
are given by restricting domains and/or codomains of bonding maps of $\G'$
if necessary.
\end{itemize}
We say that the new graph of groups is obtained from $\G$ by {\it
unkilling}.\footnote{We use the term {\it unkilling} because the
  inverse operation {\it kills a cycle}.} See \figref{f:unkilling}.

\subsection{\Splitting}\label{s:splitting}
Suppose that, for some $\bsv\in V$ and $\bse\in \hat E(\bsv)$,
$\group_{\bsv}=\group_{\bsv}'*\group_{\bsv}''$ non-trivially,
$\group_{\bse}=\group_{\bse}'*\group_{\bse}''$,
$\varphi_{\bse}(\group_{\bse}')\subset\group_{\bsv}'$,
$\varphi_{\bse}(\group_{\bse}'')\subset\group_{\bsv}''$, and for
$\bseb\in\hat E(\bsv)\setminus\{\bse\}$ either
$\varphi_{\bseb}(\group_{\bseb})\subset\group_{\bsv}'$ or
$\varphi_{\bseb}(\group_{\bseb})\subset\group_{\bsv}''$. Then, define
$\G'$ as follows:
\begin{itemize}
\item
$\bsv$ is replaced by two
vertices: $V'=V\cup\{\bsv',\bsv''\}\setminus\{\bsv\}$;
\item
the oriented edge $\{\bse,\bse^{-1}\}$ is replaced by two oriented
edges:
  $$\hat E'=\hat
  E\cup\{\bse',\bse'',(\bse')^{-1},(\bse'')^{-1}\}\setminus\{\bse,\bse^{-1}\}$$
  with $\partial_0\bse'=v'$, $\partial_0\bse''=v''$, and
  $\partial_1\bse'=\partial_1\bse''=\partial_1\bse$;
\item
for $\bseb\in \hat E(\bsv)\setminus\{\bse\}$,
$\partial_0\bseb=\bsv'$ if
$\varphi_{\bse}(\group_{\bseb})\subset\group_{\bsv}'$ and
$\partial_0\bseb=\bsv''$ if
$\varphi_{\bse}(\group_{\bseb})\subset\group_{\bsv}''$;
\item
$\group_{\bsv'}=\group_{\bsv}'$, $\group_{\bsv''}=\group_{\bsv}''$,
$\group_{\bse'}=\group_{\bse}'$, $\group_{\bse''}=\group_{\bse}''$; and
\item
bonding maps are given by restricting domains and/or codomains of
bonding maps of $\G'$ if necessary.
\end{itemize}
 We say that the new graph of groups is
obtained from $\G$ by {\it \splitting}. See \figref{f:splitting}.
\begin{figure}[ht!]\anchor{f:splitting}
\cl{\scalebox{0.75}{\input{split.pstex_t}}}
\caption{\Splitting}\label{f:splitting}
\end{figure}
Each of the operations blowing up, unpulling, unkilling, and
\splitting\ is a {\it simplification}.

\begin{proposition}
If $\G'$ is obtained from $\G$ by reducing or simplifying, then
$\pi_1(\G')$ and $\pi_1(\G)$ are isomorphic.
\end{proposition}

\begin{proof}
In the first type of reducing move,
$\pi_1(\G)\isomorphic\pi_1(\G')*_{\group_{\bse}}\group_{\bse}$ where
the map $\group_{\bse}\to\group_{\bse}$ is an isomorphism. By van
Kampen's theorem, $\pi_1(\G)\isomorphic \pi_1(\G')$. In all of the
other cases, $\G$ is obtained from $\G'$ by a Stallings fold and so
$\pi_1(\G)\isomorphic\pi_1(\G')$, see \cite[Section~2]{bf:bounding}.
\end{proof}

\begin{remark}
If $\pi_1(\G)$ is not infinite cyclic then the first type of blow up
is a composition of a second type of blow up, an unkilling, and a
reduction. Therefore, we will not have to consider the first type of
blow up.
\end{remark}

\subsection{Our case}\label{s:our case}
Given a graph of groups, we want to iteratively simplify until the
resulting graph of groups can't be simplified. In order to do this
algorithmically, we need be able to recognize when a simplification is
possible. To this end, we restrict the graphs of groups that we will
consider to the case where $\bsgraph$ is a finite graph, ie, where
$E$ is finite, and where, for $\bsv\in V$ and $\bse\in E$,
$\group_{\bsv}$ and $\group_{\bse}$ are finite rank free groups. Such
a $\G$ is a {\it finite graph of finite rank free groups}.

\subsection{Labeled graphs}\label{s:labeled graphs}
Graphs will have two uses in this paper. The first we have already
seen--these are as the underlying graphs of graphs of groups and are
denoted by $\bsgraph$'s. The other use will be to represent subgroups
of free groups and these will be denoted by $\graph$'s. We now explain
this second usage. Let $\f(\B)$ denote the free group with basis
$\B$. For $S\subset \f(\B)$, $S^{\pm 1}$ is defined to be $S\cup S^{-1}$
where $S^{-1}$ is $\{s^{-1}:s\in S\}$. A {\it labeled graph} or a {\it
$\B$--graph} is a connected graph $\graph=\graph(V,\hat
E,\op,\partial_0)$ with a {\it labeling function} $\hat E\to \B^{\pm
1}$ such that the label assigned to $\op(\e)=\e^{-1}$ is the inverse
of the label assigned to $\e$. The {\it $\B$--rose} is a $\B$--graph
$R_\B$ with one vertex and a bijective labeling function. We identify
$\pi_1(R_\B)$ with $\f(\B)$. (The homotopy class of the path formed by
the edge labeled $b$ is identified with $b$.)  There is a natural map
$\proj(\graph)\co \graph\to R_\B$ sending 1--cells to 1--cells and
preserving labels and orientations. Since $\proj(\graph)$ determines
the labeling function and {\it vice versa}, we will call
$\proj(\graph)$ the labeling function as well. The graph $\graph$ is
{\it based} if there is a distinguished vertex $*$. On the level of
fundamental groups, the image of $\proj(\graph)$ is a subgroup of
$\f(\B)$ denoted $[(\graph,*)]$. If we forget the basepoint then the
image is only defined up to conjugacy and $\graph$ determines a
conjugacy class $[[\graph]]$ of subgroups $\f(\B)$. We say that
$(\graph,*)$ {\it represents} $[(\graph,*)]$ and that $\graph$ {\it
represents} $[[\graph]]$. More generally, if $\graphseq$ is a sequence
of labeled graphs then a sequence $[[\graphseq]]$ of conjugacy classes
of subgroups of $\f(\B)$ is determined. If there are basepoints
$\vec{*}$ then a sequence $[(\graphseq,\vec{*})]$ of subgroups of
$\f(\B)$ is determined.

It is well-known, see eg \cite[Section~1.A]{ah:topology}, that a
generating set $[\graph,*]$ may be obtained as follows. Choose a
maximal tree $T$ for $\graph$ and choose orientations for the edges
not in $T$. The generating set is indexed by these oriented
edges. Specifically, the generator corresponding to the oriented edge
$\e$ is the word in $\B^{\pm 1}$ determined by reading the labels of
the loop obtained by concatenating the path in $T$ from $*$ to
$\partial_0\e$, $\e$, and the path in $T$ from $\partial_1\e$ back
to $*$.

\subsection{Stallings and Gersten representatives}\label{s:stallings}
The {\it complexity} $\complexity(\graph)$ of the labeled graph
$\graph$ is $|E(\graph)|$ and the {\it complexity}
$\complexity(\graphseq)$ of the sequence $\graphseq=\{\graph_i\}_{i\in
I}$ of labeled graphs is $\sum_{i\in I}\complexity(\graph_i)$. If
$\groupseq$ is a finite sequence of finitely generated subgroups of
$\f(\B)$, the {\it Stallings representative for $\groupseq$ with
respect to $\B$} is the sequence
$\graphseq_{\stallings}=\graph_{\stallings}(\groupseq,\B)$ of based
$\B$--graphs of minimal complexity representing $\groupseq$. We often
omit the $\B$ from the notation. If $\groupseq$ is represented by a
finite sequence $\groupseq=\{H_i=\langle S_i\rangle\}$ where each
$S_i$ is a finite set of words in $\B^{\pm 1}$, then there is an
algorithm due to Stallings \cite{st:folding} to find
$\graphseq_{\stallings}$ from $\{S_i\}$, see also \cite{fh:coherence}.

If $\conjugacyseq$ is a sequence of conjugacy classes of non-trivial
subgroups of $\f(\B)$, the {\it Stallings representative for
$\conjugacyseq$ with respect to $\B$} is the $\B$--graph
$\graph_\stallings(\conjugacyseq)$ of minimal complexity representing
$\conjugacyseq$. In fact, if $\conjugacyseq$ is represented by
$\groupseq=\{H_i=\langle S_i\rangle\}$ as above then
$\graph_\stallings(\conjugacyseq)$ is the sequence
$\core(\graph_\stallings(\groupseq))$ of cores of elements of
$\graph_{\stallings}(\groupseq)$. Recall that if $\graph$ is graph
then the {\it core} of $\graph$, denoted $\core(\graph)$, is the union
of all immersed circuits in $\graph$, see also Section~\ref{s:coring}.
In particular, there is an algorithm to find
$\core(\graph_{\stallings}(\groupseq))$ from $\{S_i\}$ as well as a
sequence $\vec{h}$ of elements of $\f(\B)$ such that
$\core(\graph_{\stallings}(\groupseq))=\graph_{\stallings}(\groupseq^{\vec{h}})$
where $\groupseq^{\vec{h}}$ is the sequence of groups obtained by
conjugating a component of $\groupseq$ with the corresponding
component of $\vec{h}$. In fact, if $\graph$ is a component of
$\graph_{\stallings}(\groupseq)$ then the corresponding component $h$
of $\vec{h}$ can be taken to be the inverse of the word read along the
shortest path from the basepoint $*$ of $\graph$ to
$\core(\graph)$. It is convenient to also allow conjugacy classes of
trivial groups in $\conjugacyseq$. Since the core of a tree is empty,
we take the Stallings representative of the conjugacy class of the
trivial group to be the empty set.

A {\it Gersten representative $\graph_{\gersten}(\conjugacyseq)$ for
$\conjugacyseq$} is a sequence of $\B$--graphs of minimal complexity
among sequences of $\B$--graphs representing $\auto\conjugacyseq$ as
$\auto$ varies over $\Aut(\f(\B))$. (If $\groupseq$ represents
$\conjugacyseq$, then $\auto\groupseq$ represents
$\auto\conjugacyseq$, where $\auto$ is applied coordinate-wise.) If
$\conjugacyseq$ is represented by $\groupseq=\{H_i=\langle
S_i\rangle\}$ as above, then there is an algorithm that produces a
$\graph_\gersten(\conjugacyseq)$ as well as an automorphism $\auto$
such that
$\core(\graph_{\stallings}(\auto\groupseq))=\graph_{\gersten}(\conjugacyseq)$,
see \cite{sg:whitehead}, \cite{sk:gersten}, and also
Section~\ref{s:gersten}.

\begin{example}
If $\B=\{a,b\}$ and if $H=\langle
aaba^{-1},ab^{-1}abba^{-1}\rangle$, then the graph $(\graph,*)$
pictured in \figref{f:gersten} represents $H$. (The open arrows
denote `$a$' and the closed `$b$'.) The graph $(\graph_{\stallings},*)$ is the
Stallings representative of $H$. For the automorphism $\auto\co \f(\B)\to
\f(\B)$ given by $a\mapsto ab^{-1}$, $b\mapsto b$, a Gersten
representative $\graph_{\gersten}$ of $H$ is the core of the Stallings
representative for $\auto H$.

\begin{figure}[ht!]\anchor{f:gersten}
\cl{\scalebox{0.75}{\input{gersten.pstex_t}}}
\nocolon\caption{}\label{f:gersten}
\end{figure}
\end{example}

\begin{notation}\label{n:vee}
If $\graphseq'=\{\graph'_i\}_{i\in I}$ and
$\graphseq''=\{\graph''_j\}_{j\in J}$ are sequences of $\B$--graphs and
if $\{0\}=I\cap J$ then $\graphseq'\vee\graphseq''=\{\graph_k\}_{k\in
I\cup J}$ is a sequence of graphs with labels in $\B$ of the following
form:
$$
\graph_k=
\begin{cases}
  \graph'_0\vee\graph''_0, &\text{if $k=0$;}\\
  \graph'_i,               &\text{if $k\in I\setminus\{0\}$; and}\\
  \graph''_j,               &\text{if $k\in J\setminus\{0\}$.}
\end{cases}
$$
\end{notation}

\begin{definition} \label{d:visible}
Let $\graphseq=\{\graph_i\}_{i\in I}$ be a sequence of $\B$--graphs.
\begin{enumerate}
\item\label{i:visible 0} If there is a non-trivial partition
$\B=\B'\sqcup\B''$ such that, for each $i\in I$, the labels of
$\graph_i$ are either all in $\B'$ or all in $\B''$, then we say that
$\graphseq$ can be {\it visibly blown up}.
\item \label{i:visible 1}
Suppose that $\lab\in\B$ appears as a label in only one element $\graph$
of $\graphseq$ and that only one oriented edge $\e_0$ of
$\graph$ has label $\lab$. Suppose further that $\e_0\subset\core(\graph)$. 
\begin{enumerate}
\item\label{i:visible unpull} If $\e_0$ does
not separate $\graph$ then we say that $\graphseq$ can be {\it
visibly unpulled}.
\item\label{i:visible unkill} If $\e_0$ does separate $\graph$
then we say that $\graph$ can be {\it visibly unkilled}.
\end{enumerate}
\item\label{i:visible split} If there is a non-trivial partition
$\B=\B'\sqcup\B''$ such that $\graphseq=\graphseq'\vee\graphseq''$ where
$\graphseq'$ is a sequence of $\B'$--graphs and $\graphseq''$ is a sequence
of $\B''$--graphs then we say that $\graphseq$ can be {\it visibly \splited}.
\end{enumerate}
If $\graphseq$ can be visibly blown up, visibly unpulled, visibly
unkilled, or visibly \splited, then we say that $\graphseq$ can be {\it
visibly simplified}.
\end{definition}

After the statement of the next proposition, we can describe the
algorithm. The proof of this proposition is almost obvious, but
requires some bookkeeping which is postponed until the appendix. This
proposition is subsumed into Proposition~\ref{p:details}.

\begin{notation}\label{n:vectors}
Let $\G$ be a graph of finite rank free groups with notation as in
Section~\ref{s:graphs of groups}.  For $\bsv\in V$, let
$\groupseq(\bsv)$ 
denote the sequence of subgroups of
$\group_{\bsv}$ represented by
$\{\varphi_{\bse}(\group_{\bse}):\bse\in\hat E(\bsv)\}$.  Also let
$\conjugacyseq(\bsv)$ denote the sequence of conjugacy classes of subgroups
of $\group_{\bsv}$ represented by $\groupseq(\bsv)$. 
\end{notation}

\begin{proposition} 
Suppose that for some $\bsv\in V$,
$\graph_{\gersten}(\conjugacyseq(\bsv))$ can be visibly
simplified. Then, $\G^\output\conjugate\G$ can be algorithmically found such
that $\G^\output$ can be simplified.
\end{proposition}

\subsection{The algorithm}\label{ss:algorithm}
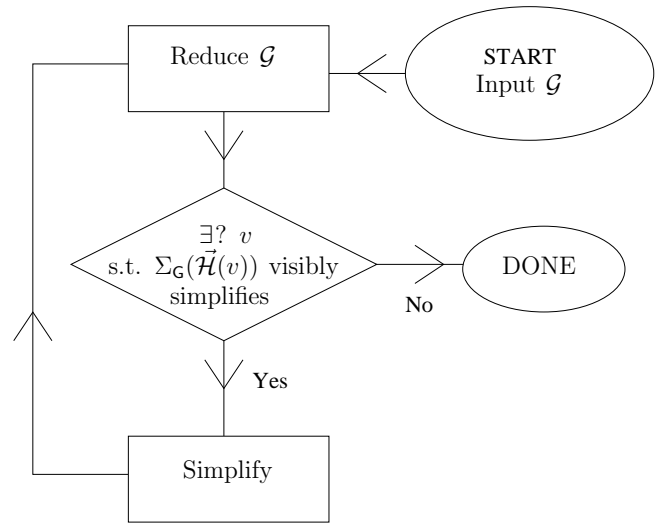
\begin{figure}[ht!]\anchor{f:flow chart}
\cl{\scalebox{0.8}{\input{flow2.pstex_t}}}
\caption{Flow chart}\label{f:flow chart}
\end{figure}

Here is the algorithm. See \figref{f:flow chart} for a flow
chart. More details on the algorithm are given in
Section~\ref{s:algorithms}, Section~\ref{s:final}, and the appendix.

\begin{description}
\item[Step 0]
Input $\G$, a finite graph of finite rank free groups.
\item[Step 1]
Reduce $\G$.
\item[Step 2] 
If, for some $\bsv\in V$,
$\graph_{\gersten}(\conjugacyseq(\bsv))$ can be visibly simplified,
then replace $\G$ by a simplified conjugate and return to Step
1. Else, done.
\end{description}

The main result of the paper is:

\begin{theorem}\label{t:main 2}
Suppose that a finite graph $\G$ of finite rank free groups is input
into the above algorithm and that $\G^\output$ is output. Then,
the decomposition of $\pi_1(\G)$ determined by the edges
of $\G^\output$ with trivial stabilizer is the Grushko decomposition
of $\pi_1(\G)$.
\end{theorem}

The proof of Theorem~\ref{t:main 2} is found in Section~\ref{s:final}.

\subsection{Relative version}\label{ss:relative}
In this section we describe a relative version of Theorem~\ref{t:main
2}. If $H$ is a subgroup of a group $\group$, then we say that
$\group$ is {\it freely decomposable rel $H$} if there is a free
decomposition $\group=\group'*\group''$ with $H\subset\group'$ and
$\group''$ non-trivial. Otherwise, $\group$ is {\it freely
indecomposable rel $H$}. The relative version of Grushko's theorem
(Theorem~\ref{t:grushko}) is:

\begin{theorem}
Suppose $H$ is a subgroup of the finitely generated group
$\group$. Then, $\group$ is a free product
$\group_H*\group_1*\dots*\group_n*F$ where $H\subset \group_H$,
$\group_H$ is freely indecomposable rel $H$, $\group_i$ for $1\le i\le
n$ is freely indecomposable and not free, and $F$ is a finite rank
free group.
\end{theorem}

The subgroup $\group_H$ is unique in this {\it the Grushko
decomposition of $\group$ rel $H$.}  Up to reordering and conjugation,
the $\group_i$, $1\le i\le n$, are unique. Also, the rank of $F$ is an
invariant of the pair $H\subset \group$.

Suppose now that $\G$ is a finite graph of finite rank free groups,
that $\bsv_0\in V$ has valence one with incident edge $\bse_0$, and
that $\varphi_{\bse_0}$ is an isomorphism.  We are going to describe a
slight modification of the algorithm of Section~\ref{ss:algorithm}
that produces the Grushko decomposition of $\pi_1(\G)$ rel
$\group_{\bsv_0}$. Intuitively, in the modified algorithm we only
reduce or visibly simplify only if the special edge group
$\group_{\bse_0}$ is unchanged. Specifically, we modify the algorithm
as follows.

Steps~0 and 1 are replaced by:
\begin{description}
\item[Step 0$'$] Input $\G$, a finite graph of finite rank free
groups as above.
\item[Step 1$'$]
Reduce $\G$ rel $\bsv_0$, ie, apply the reducing moves displayed in
\figref{f:reduce} only if $\bse\notin\{\bse_0,\bse_0^{-1}\}$.
\end{description}

To describe the modification of Step~2, we need a definition. In
Definition~\ref{d:visible}(\ref{i:visible 1}), the component $\graph$
of $\graphseq$ is {\it special}. In
Definition~\ref{d:visible}(\ref{i:visible split}), the component of
$\graphseq$ corresponding to $\graph_0=\graph_0'\vee\graph_0''$ in
Notation~\ref{n:vee} is {\it special}. In
Definition~\ref{d:visible}(\ref{i:visible 0}), none of the components
of $\graphseq$ are special. The components of
$\graph_{\gersten}(\conjugacyseq(\bsv))$ are parametrized by the set
$\hat E(\bsv)$ of edges incident to $\bsv$. If
$\graph_{\gersten}(\conjugacyseq(\bsv))$ can be visibly simplified,
then the edge corresponding to the special component is {\it
special}. If the special edge $\bse$ is not in
$\{\bse_0,\bse_0^{-1}\}$, then the resulting simplification will not
change $\group_{\bse_0}$. Step~2 of the algorithm is replaced by:
\begin{description}
\item[Step~2\,$'$] If, for some $\bsv\in V$,
$\graph_{\gersten}(\conjugacyseq(\bsv))$ can be visibly simplified and
the special edge is not in $\{\bse_0,\bse_0^{-1}\}$, then replace $\G$
by a simplified conjugate and return to Step~$1'$. Else, done.
\end{description}

\begin{theorem}\label{t:relative}
Suppose that $\G$ as above is input into the modified algorithm and
that $\G^\output$ is output. Then, the decomposition of $\pi_1(\G)$
determined by the edges of $\G^\output$ with trivial stabilizer gives
the Grushko decomposition of $\pi_1(\G)$ rel $\group_{\bsv_0}$.
\end{theorem}
\noindent
The proof of the relative version requires only minor notational changes
to the proof of Theorem~\ref{t:main 2} and is left to the reader.

\section{Laminated square complexes and models}\label{s:square complexes}
This section contains a discussion of certain laminated two complexes
called models whose 2--cells are squares. For an interesting study of
complexes built from squares see \cite{bw:squares}.

Let $I$ denote the unit interval $[0,1]$. An {\it $n$--cube} is a
metric space isometric to $I^n$.  A metric space is a {\it cube} if it
is an $n$--cube for some $n$. A {\it cube complex} is a union of cubes
glued by isometries of faces. A finite dimensional cube complex $X$
admits a maximal metric such that the inclusion $C\to X$ is a local
isometry for each cube $C$ of $X$ \cite{bh:nonpositive}. A {\it
square} is a 2--cube, and we only have need to consider {\it square
complexes}, ie, cube complexes of dimension at most two. A {\it
graph} is a 1--dimensional cube complex\footnote{1--dimensional
$CW$--complexes and 1--dimensional cube complexes are both called {\it
graphs}. Since we will only be using combinatorial properties, the
distinction is not important to us.}, a {\it tree} is a simply
connected graph and a {\it forest} is a disjoint union of trees.

A decomposition of a square is {\it standard} if it is induced by
projection to a codimension--1 face. A decomposition of a 1--cube is
{\it standard} if all decomposition elements are points. It is {\it
trivial} if the only decomposition element is the 1--cube itself. A {\it
laminated square complex} is a simply connected square complex $X$
with a decomposition $\D$ such that:
\begin{description}
\item[\rm(M1)]
The link of every vertex of $X$ is a flag complex.
\item[\rm(M2)] For each square $C$ of $X$, the induced decomposition of
  $C$ is standard. In other words, the decomposition of $C$ whose
  elements are the components of $C$ intersected with elements of $\D$
  is standard. For each 1--cube of $X$, the induced decomposition is
  either standard or trivial.
\end{description}
In this context, (M1) means that every link is a simplicial graph
with no circuits of length three. A decomposition element is also
called a {\it leaf}.

\begin{proposition}
Let $(X,\D)$ be a laminated square complex. Then,
\begin{enumerate}
\item
$X$ is contractible and
\item
leaves are forests.
\end{enumerate}
\end{proposition}

\begin{proof}
(M1) implies that the the metric on $X$ is $CAT(0)$ \cite{bh:nonpositive}.
Hence, (1). A vertex of a square with a standard decomposition is
contained in exactly one edge that is contained in a leaf. So, links
of vertices in $X$ are bipartite and leaves are totally geodesic. In
particular, leaves are 1--dimensional and contractible.
\end{proof}

A {\it model} is a laminated square complex $(X,\D)$ such that:
\begin{description}
\item[\rm(M3)]
Leaves are connected.
\end{description}

The next proposition is immediate from definitions.

\begin{proposition}
Let $(X,\D)$ be a laminated square complex and let $\induce(\D)$ be the
decomposition of $X$ whose elements are the connected
components of elements of $\D$. 
Then, $(X,\induce(\D))$ is a model.\qed
\end{proposition}

\begin{definition} \label{d:induced}
We say that $(X,\induce(\D))$ is {\it induced by} $(X,\D)$.
\end{definition}

\begin{proposition}
If $(X,\D)$ is a model, then the decomposition space $X/\D$ is a tree.
\end{proposition}

\begin{proof}
Since a leaf is totally geodesic, it intersects each square $C$ in a
connected set. In particular, the decomposition $C\cap\D$ of $C$
obtained by intersecting $C$ with elements of $\D$ is standard,
$I= C/C\cap\D\to X/\D$ is injective, and $X/\D$ is naturally a
graph. Leaves are connected, and so $X\to X/\D$ is
$\pi_1$--surjective. Thus, $X/\D$ is a tree.
\end{proof}

If $(X,\D)$ is a model, we say that {\it $X$ is a model for the tree
$S=X/\D$} or that {\it $X\to S$ is a model}. The preimage in $X$ of
$s\in S$ is $X_s$.  If $\e$ is an edge of $S$ and if
$\s,\s'\in\mathring{\e}$, then $X_\s$ and $X_{\s'}$ have the same
isomorphism type $X_\e$.  We will sometimes abuse notation and
identify $X_\e$ with $X_s$ for $s\in\mathring{\e}$. 

\begin{notation}
If the group $G$ acts on the set $X$ and if $S$ is a subset of $X$,
then $G_{X,S}$ is the stabilizer of $S$, ie, the subgroup elements
$g\in G$ such that $g(S)=S$. If $S=\{s\}$ then we also write $G_{X,s}$ for
$G_{X,S}$. We will suppress the $X$ if the space is understood.
\end{notation}

An action of a group $\group$ on $(X,\D)$ is an (isometric) action of
$\group$ on $X$ permuting cubes and decomposition elements. In this
case, $(X,\D)$ is a {\it $G$--model}. There is an induced action of
$\group$ on $S=X/\D$. The quotients $X/\group$ and $S/\group$ are
denoted $\bar X$ and $\bar S$ respectively.  We say that {\it $G$ acts
without inversions} if, for all cubes $C\subset X$, $G_C$ fixes $C$
pointwise. By subdividing $X$, we may arrange that $G$ acts without
inversions. Hence, {\sl we always assume that our actions are without
inversions.}  Note that the space $X_\s$ is $\group_s$--invariant and
$X_\e$ is $\group_\e$--invariant.

\section{Trees}\label{s:trees}
We review some tree basics. A $\group$--tree $S$ is {\it minimal} if it
has no proper $\group$--invariant subtrees. It is {\it trivial} if has
a fixed point. In this case we also say that the $\group$--action is
{\it elliptic.} If $S$ has a unique minimal invariant
$\group$--subtree, then it is denoted $S_\group$. This occurs, for
example, if $\group$ contains a {\it hyperbolic} element, that is an
element that fixes no point of $S$ \cite{hb:remarks}. A {\it morphism
$S\to T$} of $\group$--trees is a simplicial $\group$--map. It is {\it
\proper}\ if no edge is mapped to a point. If there is a morphism $S\to
T$, then {\it $S$ resolves $T$}.

If $\edge$ is an edge of $S$, then we say that {\it $S'$ is obtained
from $S$ by collapsing $\edge$} if $S'$ is the result of equivariantly
collapsing $\edge$. A morphism $S\to S'$ is a {\it collapse} if $S'$
is obtained from $S$ by iteratively collapsing edges.  If edges
$\edge_1$ and $\edge_2$ of $S$ share the vertex $\s$, then we say that
{\it $S'$ is obtained from folding $\edge_1$ and $\edge_2$} if $S'$ is
the result of equivariantly identifying $\edge_1$ and $\edge_2$ with
an isometry fixing $\s$. The resulting morphism $S\to S'$ is a {\it
fold}.

\section{Examples of Models}\label{s:examples}
\begin{example} 
First a non-example. Glue two squares along three sides and laminate
so that restrictions to squares are standard and so that the unglued
sides form a leaf. The result is not a laminated square complex even
though it is contractible (there are vertices whose links consist of
distinct edges with the same endpoints--such a link is not a
simplicial graph). Notice that not all leaves are trees.
\end{example}

\begin{example}
The quotient $\bar X\to\bar S$ of a 
$\group$--model $X\to S$ with
$\group$ free of rank 2 is depicted in \figref{f:model}. The
preimage in $\bar X$ of a point in the interior of the edge of
$\bar S$ is isomorphic to a circle. The preimage of the vertex is
a `pair of eyeglasses'. The stabilizer of a vertex of $S$ is free of
rank two; the stabilizer of an edge of $S$ is infinite cyclic.

\begin{figure}[ht!]\anchor{f:model}
\cl{\scalebox{0.6}{\input{model.pstex_t}}}
\nocolon\caption{}
\label{f:model}
\end{figure}
\end{example}

\begin{example}
If $S$ and $T$ are $\group$--trees, then 
$$S\times T=\cup\{e\times f : e\mbox{ is an edge of }S,\ f\mbox{ is
an edge of }T\}$$ is a union of squares with projection maps
$\treeproj(S)\co S\times T\to S$ and $\treeproj(T):S\times T\to T$.  The
induced decomposition with quotient $S$ (respectively $T$) gives
$S\times T$ the structure of a model for $S$ (respectively $T$).
\end{example}

\begin{example}\label{e:subcomplex}
Let $(X,\D)$ be a $G$--model and let $Y$ be a simply connected
$\group$--subcomplex of $X$. Then, $Y$ with the decomposition
$\D(Y)=\{D\cap Y\mid D\in\D\}$ is a laminated square complex and
$Y/\D(Y)\to X/\D$ is an inclusion. If we let $\induce(\D)(Y)$ be the
decomposition of $Y$ induced by $\D(Y)$ (see
Definition~\ref{d:induced}), then $(Y,\induce(\D)(Y)))$ is a model and
$Y/\induce(\D)(Y)\to X/\D$ is a morphism. We call
$\induce(\D)(Y)$ the {\it restricted decomposition of $Y$}.
\end{example}

\begin{example}[Main Example]\label{e:main}
Suppose that $S$ and $T$ are $\group$--trees and that, for each $\s\in
S$, there is a unique minimal $\group_\s$--invariant subtree $T_\s$ of
$T$. The union $X_S(T)=\bigsqcup_{\s\in S}(\{s\}\times T_\s)$ is a
subcomplex of $S\times T$ and is simply-connected (being a union of
simply-connected spaces along simply-connected spaces). $X_S(T)\to S$
is a model as is $X_S(T)\to\induce(T)$ where $\induce(T)$ is the
quotient of the decomposition induced from $X_S(T)\subset S\times T\to
T$ by restriction (see Example~\ref{e:subcomplex}). If $\overline
S=S/G$ and all $\overline T_s=T_s/G_s$ are compact, then
$\overline{X_S(T)}=X_S(T)/G$ is also compact.

\begin{proposition}\label{p:minimal}
$\induce(T)$ as in Example~\ref{e:main} is minimal.
\end{proposition}

\begin{proof}
We may identify $T_s$ with the image of the injection $\{s\}\times
T_s\to \induce(T)$. 
By construction,
$\induce(T)=\cup_{s\in S}T_s$. If $\induce(t)\in\induce(T)$ is not
contained in an invariant $G$--subtree $R$ of $\induce(T)$ and if
$\induce(t)\in T_s$, then $R\cap T_s$ is a proper $G_s$--invariant
subtree of $T_s$, contradiction.
\end{proof}
\end{example}

\section{Operations on models}\label{s:operations on models}
{\sl In this section, we assume that $\group$ is a group, $S$ is a
$\group$--tree, and $X\to S$ is a model.} We will describe operations on
$X$. In each case, the result $\result(X)$ is a model for $\result(S)$
where $\result(S)$ resolves $S$. The operations are geometric
generalizations of the simplifications of Section~\ref{s:algorithm}.

\subsection{0--Simplifying}
Let $C=I\subset X$ be a 1--cube and set $C_0=\{0\}$. Suppose that 
\begin{itemize}
\item
$C$ meets cubes of $X$ other than faces of $C$ only in $\{1\}$; and
\item
the restriction of the decomposition to $C$ is $\{C\}$.
\end{itemize}
Let $\result(X)$ be the result of equivariantly replacing $C$ by
$\{1\}$, ie, $$\result(X)=X\setminus\cup_{g\in G}g\cdot [0,1).$$ We
say that $\result(X)$ with the restricted decomposition is obtained
from $X$ by {\it 0--simplifying $C$ from $C_0$}. Here
$\result(S)=S$. See \figref{f:0collapse}.

\begin{figure}[ht!]\anchor{f:0collapse}
\cl{\scalebox{0.65}{\input{0collapse.pstex_t}}}
\nocolon\caption{}
\label{f:0collapse}
\end{figure}

\subsection{I--Simplifying}
Let $C=I\times I^n\subset X$ ($n=0$ or 1) be a cube and set
$C_0=\{0\}\times I^n$. Suppose that 
\begin{itemize}
\item $C$ meets cubes of $X$ other than
faces of $C$ only in $\{1\}\times I^n$; and 
\item
$C_0$ is a decomposition
element. 
\end{itemize}
Let $\result(X)$ be the result of equivariantly replacing
$C$ by $\{1\}\times I^n$.
We say that $\result(X)$ with the restricted decomposition is obtained
from $X$ by {\it I--simplifying $C$ from $C_0$}. $S\to\result(S)$
is a collapse. See \figref{f:1collapse}.
\begin{figure}[ht!]\anchor{f:1collapse}
\cl{\scalebox{0.75}{\input{1collapse.pstex_t}}}
\nocolon\caption{}
\label{f:1collapse}
\end{figure}

\subsection{II--Simplifying}
Let $C=I^2\subset X$ be a square and set $C_0=\{0\}\times I$. Suppose
that 
\begin{itemize}
\item
$C$ meets cubes of $X$ other than faces of $C$ only in
$$(\{1\}\times I)\cup (I\times\{ 1\});$$ 
\item
it is not possible to I--simplify from $C_0$, ie, $[0,1)\times\{1\}$
meets a cube other than a face of $C$; and
\item
$C_0$ is an element of the decomposition restricted to $C$.
\end{itemize}
The model
$\result(X)$ with restricted decomposition elements obtained by
equivariantly replacing $C$ by
$$(\{1\}\times I)\cup (I\times\{ 1\})$$ is the result of {\it
II--simplifying $C$ from $C_0$}. Note that $\result(S)=S$. See
\figref{f:2collapse}.
\begin{figure}[ht!]\anchor{f:2collapse}
\cl{\scalebox{0.63}{\input{2collapse.pstex_t}}}
\nocolon\caption{}
\label{f:2collapse}
\end{figure}

\subsection{III--Simplifying}
Let $C=I^2\subset X$ be a square and set $C_0=\{0\}\times I$. Suppose
that 
\begin{itemize}
\item
$C$ meets cubes of $X$ other than faces of $C$ only in
$(\partial{I}\times I)\cup (I\times\{ 1\});$ 
\item
it is not possible to I--, or II--simplify from $C_0$; and
\item
$C_0$ is an element of the decomposition restricted to $C$.
\end{itemize}
The model $\result(X)$ with restricted decomposition
elements obtained by equivariantly replacing $C$ by
$(\partial{I}\times I)\cup (I\times\{ 1\})$ is the result of {\it
III--simplifying $C$ from $C_0$.}  Note that $\result(S)\to S$ is a
non-trivial fold. See \figref{f:3collapse}.
\begin{figure}[ht!]\anchor{f:3collapse}
\cl{\scalebox{0.63}{\input{3collapse.pstex_t}}}
\nocolon\caption{}
\label{f:3collapse}
\end{figure}

\subsection{Blowing up}
If there is a cube $C=I\subset X$ such that $\mathring{I}$ meets no
cube other than faces of $C$ and such that the decomposition
restricted to $C$ is $\{C\}$ then we may refine the decomposition by
$\group$--equivariantly replacing the decomposition element $X_C$
containing $C$ with
$$\{\{g\cdot t\}\mid t\in\mathring{I},\ g\in
G_{X_C}\}\sqcup\{\mbox{components of }X_{C}\setminus
G_{X_C}\cdot\mathring{I}\}.$$ We say that $\result(X)\to\result(S)$ is
obtained from $X\to S$ by {\it blowing up $C$}. The induced map
$\result(S)\to S$ collapses to points the edges of $\result(S)$
corresponding to the orbit of $C$, explaining the term ``blowing
up''. See \figref{f:blowup}.
\begin{figure}[ht!]\anchor{f:blowup}
\cl{\scalebox{0.63}{\input{blowup.pstex_t}}}
\nocolon\caption{}
\label{f:blowup}
\end{figure}

\begin{definition}\label{d:simplifying trees}
Let $\result(S)\to S$ a morphism. If there is a model $Y\to S$ and a
cube $C$ of $Y$ with face $C_0$ such that 0--simplifying $C$ from $C_0$
yields $\result(Y)\to\result(S)$, then we say that $\result(S)$ is
obtained from $S$ by {\it 0--$\groupo$--simplifying} or equivalently by
{\it 0--simplifying over $\groupo$} where $\groupo$ is the stabilizer
of $C_0$ in $Y$. The definitions for
I--$\groupo$--, II--$\groupo$--, and III--$\groupo$--simplifying are analogous.
If there is a model $Y\to S$ with a cube $C$ such that
$\result(Y)\to\result(S)$ is the result of blowing up $C$, then we say
that $\result(S)$ is obtained from $S$ by {\it $\groupo$--blowing up}
or equivalently by {\it blowing up over $\groupo$} where $\groupo$ is
the stabilizer of $C$ in $Y$. In each of these cases, we say that
$\result(S)$ is obtained from $S$ by {\it simplifying over $\groupo$}
or just by {\it simplifying} if $\groupo$ is understood.
\end{definition}

\begin{remark}
Since the identity map $S\to S$ is an example of a model, blowing up
as defined in Section~\ref{s:blowing up} is an example of
{$\trivialgroup$}--blowing up.
\end{remark}

\section{Generalized Shenitzer--Swarup}\label{s:gen s-s 1}
In the next lemma, we use the notation of Example~\ref{e:main}.
\begin{lemma}\label{l:s-s for X_S(T)}
Suppose that $\group$ is a group and that $S$ and $T$ are
$\group$--trees such that
\begin{enumerate}
\item
$\overline S$ is compact;
\item
for each $\s\in S$ there is a unique minimal cocompact $\group_s$--subtree
$T_\s$ of $T$; and
\item
for each edge $\et\subset T$, the action of $\group_\et$ on $S$ is elliptic.
\end{enumerate}
Then, there is a sequence
$$\{X_s(T)=X_0\to S=S_0, X_1\to S_1, \cdots, X_N\to S_N\}$$ of I--,
II--, and III--simplifications such that $(X_N)_{\induce(\et)}$ is a
point for each edge $\induce(\et)$ of $\induce(T)$. Further, all the
simplifications are over subgroups of edge stabilizers
of $T$.
\end{lemma}

\begin{proof}
Recall from Example~\ref{e:main} that $X_0\to \induce(T)$ is obtained
by restricting $S\times T\to T$.  It is not possible to 0--simplify
$X_0\to S_0$. Indeed, in order to 0--simplify $X_0$ there would have to
be a 1--cube as in the definition of 0--simplifying. By the
construction of $X_0$, the restriction of the decomposition giving
$X_0\to\induce(T)$ to this 1--cube is standard. This is impossible
since, by Proposition~\ref{p:minimal}, $\induce(T)$ is minimal.
Further, if $X_i$ is obtained from $X_0$ by a sequence of I--, II--, and
III--simplifications, then the restriction to $X_i$ of the
decomposition giving $X_0\to\induce(T)$ still has decomposition space
$\induce(T)$. (It is the decomposition space $S_i$ of the restriction
of $S\times T\to S$ that can change.) In particular, it is also not
possible to 0--simplify $X_i$.

Suppose we have constructed the sequence $$\{X_S(T)=X_0\to S=S_0,\
X_1\to S_1,\ \cdots,\ X_i\to S_i\}.$$ We will describe how to proceed.
Let $\induce(\t)$ be a point in the interior of an edge of
$\induce(T)$. The preimage $(X_i)_{\induce(\t)}$ of $\induce(\t)$
under $X_i\to\induce(T)$ is a $\group_{\induce(\t)}$--subtree of
$S\times\{\induce(t)\}$. By (1) and (2), $X_0$, and so also $X_i$, is
cocompact.  Therefore, $(X_i)_{\induce(t)}/\group_{\induce(t)}$ is
compact. Since $T'$ resolves $T$, by (3) the action of
$\group_\induce(t)$ on $S\times\{\induce(t)\}$, and hence also on
$X_\induce(t)$ is elliptic. We see that
$X_\induce(t)/\group_\induce(t)$ is a finite tree. If this finite tree
is not a single vertex then $X_\induce(t)$ contains a valence one
vertex whose
stabilizer equals the stabilizer of the incident edge.

If there is such a valence one vertex, then this vertex is contained
in a cube $C^i_0=I\subset X_i$ that projects to a point in $S_i$. In
this case, simplify from $C^i_0$ to obtain $X_{i+1}\to S_{i+1}$. Stop
if, for each edge $\induce(\et)$ of $\induce(T)$,
$\induce(T)_{\induce(\et)}$ is a point. 

The process must eventually stop since there are only finitely many
$\group$--orbits of cubes in $X_S(T)$.

The final claim of the lemma follows from the observation that
$C^i_0$ projects to an edge of $\induce(T)$ and so the stabilizer
$\group^i_0$ of $C^i_0$ fixes an edge of $\induce(T)$. Since
$\induce(T)$ resolves $T$, $\group^i_0$ fixes an edge of $T$ as well.
\end{proof}

\begin{theorem}\label{t:gen s-s 1}
Let $S$ be a cocompact $\group$--tree with finitely generated edge
stabilizers and with $\group$ finitely generated. Suppose that
$\group$ splits over a finite group. Then, $S$ may be iteratively I--,
II--, and III--simplified and then blown up to a $\group$--tree
$\result(S)$ such that the decomposition of $\group$ given by edges of
$\result(S)$ with finite stabilizer is non-trivial. Further, the
simplifications and blow ups are all over finite groups. In
particular, all point stabilizers of $\result(S)$ are finitely
generated.
\end{theorem}

\begin{proof}
Choose $T$ to be a minimal $\group$--tree with one orbit of edges and
with finite edge stabilizers. If an edge stabilizer of $S$ is finite
then we may set $\result(S)=S$ and we are done. We may assume then
that the edges stabilizers of $S$ are infinite.

Since edge stabilizers of $S$ are finitely generated and since
$\group$ is finitely generated, for each $s\in S$, $\group_s$ is
finitely generated, see for example \cite[Lemma~32]{dc:book}. The edge
stabilizers of $T$ are finite and by assumption $\group_s$ is infinite
and so either $\group_s$ is contained in a unique vertex stabilizer of
$T$ or some element of $\group_s$ acts hyperbolically on $T$. In
particular, there is a unique minimal cocompact $\group_s$--subtree
$T_s$ of $T$. Therefore we may apply Lemma~\ref{l:s-s for X_S(T)} to
simplify $X_S(T)$ to obtain $X_N\to S_N$.

Blow up $X_N\to S_N$ to obtain $\result(X)\to \result(S)$. Since
$(X_N)_\induce(f)$ is a point for each edge $\induce(\et)$ of
$\induce(T)$, $\result(S)$ resolves $\induce(T)$. By
Proposition~\ref{p:minimal}, $\induce(T)$ is minimal and so $\result(S)$ is
non-trivial.
\end{proof}

\begin{corollary}[Generalized Shenitzer--Swarup]\label{c:gen s-s}
Let $S$ be a minimal $\group$--tree with finitely generated edge
stabilizers.  Suppose that $\group$ splits over $\trivialgroup$. Then,
$S$ may be iteratively $\trivialgroup$--simplified to a tree
$\result(S)$ such that the decomposition of $\group$ determined by the
edges of $\result(S)$ with trivial stabilizer is non-trivial.\qed
\end{corollary}

The focus of this paper is on splittings over $\trivialgroup$, ie,
on free decompositions. In a future paper, we plan to explore
splittings over small groups. Here is a sample analogue of
Lemma~\ref{l:s-s for X_S(T)} in that setting. Again, we use the
notation of Example~\ref{e:main}.

\begin{theorem}
Suppose that $\group$ is a freely indecomposable 
group. Suppose that $S$ and $T$ are $\group$--trees such that
\begin{enumerate}
\item
$\overline S$ is compact;
\item
for each $s\in S$, there is a unique minimal cocompact
$\group_s$--subtree $T_s$ of $T$; and
\item
edge stabilizers of $T$ are infinite cyclic and $T$ has one orbit of
edges.
\end{enumerate}
Then, $X_S(T)\to S$ may be iteratively simplified to
$\result(X)\to \result(S)$ where, for each edge $\induce(\et)$ of
$\induce(T)$, $\result(X_{\induce(\et)})$ is either a point with
infinite cyclic stabilizer or a line with infinite cyclic stabilizer
with generator acting by a non-trivial translation. Further, these
simplifications are over $\trivialgroup$ or $\Z$.
\end{theorem}

\begin{proof}
Since $\induce(T)$ resolves $T$, $\induce(T)$ is minimal, and $\group$
is freely indecomposable, it follows that the edge stabilizers of
$\induce(T)$ are infinite cyclic. Thus, for $\induce(\et)$ an edge of
$\induce(T)$, $(X_S(T))_\induce(\et)$ is a $\Z$--tree. If this tree is
not a point or a line, then it has a valence one vertex and a
simplification is possible. Iterate.
\end{proof}

\section{Algebraic consequences}\label{s:algebra}
This section will be needed for algorithmic questions. We use the
notation of Section~\ref{s:operations on models}. The goal is to
describe the effect of simplifying on edge and vertex stabilizers.

\begin{definition}
Let $S$ be a $\group$--tree and let $X\to S$ be a model with a cube $C$
with face $C_0$ such that $\result(X)\to\result(S)$ is the result of
III--simplifying $C$ from $C_0$.  Further, let $\es$ be the image in
$S$ of $C$, let $\s_0$ be the image of $C_0$, and let $\s_\es\in
\mathring{\es}$. Set $C_{\s_\es}=X_{\s_\es}\cap C$. Denote by
$\overline C_{\s_\es}$ the image of $C_{\s_\es}$ in $\overline X_{\s_\es}$ and
by $\overline C_0$ the image of $C_0$ in $\overline X_{\s_0}$. There are
three cases.
\begin{enumerate}
\item\label{i:both} $\overline C_{\s_\es}$ separates $\overline X_{\s_\es}$
and $\overline C_0$ separates $\overline X_{\s_0}$.
\item\label{i:one} $\overline C_{\s_\es}$ separates $\overline
X_{\s_\es}$, but $\overline C_0$ does not separate $\overline X_{\s_0}$.
\item\label{i:neither} $\overline C_{\s_\es}$ does not separate
$\overline X_{\s_\es}$ and $\overline C_0$ does not separate
$\overline X_{\s_0}$.
\end{enumerate}
Let $\groupo$ be the stabilizer in $X$ of $C_0$.
In Case~(\ref{i:both}), we say the simplification is a {\it
$\groupo$--\splitting}, in Case~(\ref{i:one}) a {\it $\groupo$--unkilling}, and
in Case~(\ref{i:neither}) a {\it $\groupo$--unpulling}.
\end{definition}

For the moment, we forget models and make some purely algebraic
definitions. Here $\group_0$ is a subgroup of
the group $\group$ and $\conjugacyseq$ is a sequence of conjugacy
classes of subgroups of $\group$.

\begin{definition}\label{d:algebraic splitting}
If there are subgroups $\groupi\subset\group$ and
$\groupii\subset\group$ containing $\groupo$ such that
\begin{itemize}
\item
$\group=\groupi*_{\groupo} \groupii$, ie, the natural map
$\groupi*_{\groupo} \groupii\to\group$ is an isomorphism;
\item
some $\conjugacyo\in\conjugacyseq$ has a representative
$\subgroupo\in\conjugacyo$ with subgroups $\subgroupoi\subset\subgroupo$ and
$\subgroupoii\subset\subgroupo$ satisfying
\begin{itemize}
\item
$\groupo\subset\subgroupoi\subset\groupi$;
\item
$\groupo\subset\subgroupoii\subset\groupii$; and
\item
$\subgroupo=\subgroupoi*_{\groupo} \subgroupoii$; and
\end{itemize}
\item
for all $\conjugacy\not=\conjugacyo$ in $\conjugacyseq$, there is
$\subgroup\in\conjugacy$ such that either $\subgroup\subset\groupi$
or $\subgroup\subset\groupii$
\end{itemize}
then we say that {\it $\conjugacyseq$ can be $\groupo$--\splited\ in $\group$.}
\end{definition}

\begin{definition}\label{d:algebraic unkilling}
If there is a subgroup $\groupi\subset\group$ containing $\groupo$, a
monomorphism $h\co \groupo\to\groupi$, and $t\in\group$ such that
\begin{itemize}
\item
$\group=\groupi*_h=\langle \groupi,t\mid tgt^{-1}=h(g),\ g\in \groupo\rangle$;
\item
some $\conjugacyo\in\conjugacyseq$ has a representative
$\subgroupo\in\conjugacyo$ with subgroups $\subgroupoi$ and
$\subgroupoii$ satisfying
\begin{itemize}
\item
$\groupo\subset\subgroupoi\subset\groupi$;
\item
$h(\groupo)\subset\subgroupoii\subset\groupi$; and
\item
$\subgroupo=\subgroupoi*_{\groupo} t^{-1}\subgroupoii t$; and
\end{itemize}
\item
for all $\conjugacy\not=\conjugacyo$ in $\conjugacyseq$ there is $\subgroup\in\conjugacy$ with $\subgroup\subset\groupi$
\end{itemize}
then we say that {\it $\conjugacyseq$ can be $\groupo$--unkilled in $\group$.}
\end{definition}

\begin{definition}\label{d:algebraic unpulling}
If there is a subgroup $\groupi\subset\group$ containing $\groupo$, a
monomorphism $h\co \groupo\to\groupi$, and $t\in\group$ such that
\begin{itemize}
\item
$\group=\groupi*_h$;
\item
some $\conjugacyo\in\conjugacyseq$ has a
representative $\subgroupo\in\conjugacyo$ with a subgroup
$\subgroupoi$ satisfying
\begin{itemize}
\item
$t\in\subgroupo$;
\item
$\subgroupoi\subset\groupi$;
\item
$\groupo\subset\subgroupoi$;
\item
$h(\groupo)\subset\subgroupoi$; and
\item
$\subgroupo=\subgroupoi*_{h}$; and
\end{itemize}
\item
for all $\conjugacy\not=\conjugacyo$, there is $\subgroup\in\conjugacy$
such that $\subgroup\subset\groupi$
\end{itemize}
then we say that {\it $\conjugacyseq$ can be $\groupo$--unpulled in
$\group$.}

If $\conjugacyseq$ can be $\groupo$--\splited, $\groupo$--unkilled, or
$\groupo$--unpulled then we say it can be {\it III--$\groupo$--simplified}.
\end{definition}

Recall that if $\s$ is a vertex of the $\group$--tree $S$ then
$\conjugacyseq(s)$ denotes the sequence of conjugacy classes of
subgroups of the stabilizer of $\s$ represented by the stabilizers of
edges incident to $\s$. The sequence is indexed by the oriented edges
$\overline S=S/\group$ that are incident to the image in $\overline S$
of $\s$.

\begin{lemma}\label{l:simplify characterization}
Let $S$ be a $\group$--tree. If $S$ can be $\groupo$--\splited,
$\groupo$--unkilled, or $\groupo$--unpulled then there is a vertex
$\s$ of $S$ such that $\conjugacyseq(\s)$ can be $\groupo$--\splited,
$\groupo$--unkilled, or $\groupo$--unpulled in $\group_\s$.
\end{lemma}

\begin{proof}
Assume that $S$ can be III--$\groupo$--simplified. Let $X\to S$ be a
model with a cube $C=I\times I$ and face $C_0=\{0\}\times I$ such that
III--simplifying $C$ from $C_0$ produces $\result(X)\to\result(S)$. In
particular, $\groupo$ is the stabilizer of $C_0$.  Let $\s_0$ be the
image of $C_0$ in $S$, let $\es$ be the image of $C$ in $S$, and let
$\subgroupo$ be the stabilizer of $\es$. Choose
$\s_{\es}\in\mathring{\es}$, and set $C_{\s_\es}=X_{\s_\es}\cap
C$.

The desired splitting of $\group_{\s_0}$ is obtained by collapsing all
edges of the $\group_{\s_0}$--tree $X_{\s_0}$ that are not in the orbit
of $C_0$. The desired $\subgroupo$--tree is obtained by
collapsing all edges of the $\group_0$--tree $X_{\s_\es}$ that are not in
the orbit of $C_{\s_\es}$. Thus, $\conjugacyseq(\s)$ can be
III--$\groupo$--simplified.
\end{proof}

\begin{definition}
Suppose that $\group$ is a group and that $\conjugacyseq$ is a
sequence of conjugacy classes of subgroups of $\group$. Suppose that
$\group=\groupi*_{\groupo} \groupii$ or
$\group=\groupi*_{\groupo}$ and that, for all
$\conjugacy\in\conjugacyseq$, $\conjugacy$ is conjugate into
$\groupi$ or $\groupii$. Then we say that $\conjugacyseq$ can be {\it
$\groupo$--blown up in $\group$.}
\end{definition}

The proofs of the Lemmas~\ref{l:blowing up} and \ref{l:I-simplifying}
are very similar to that of Lemma~\ref{l:simplify characterization}
and are not provided.
\begin{lemma}\label{l:blowing up}
Let $\group$ be a group and let $S$ be a $\group$--tree. If $S$ can be
$\groupo$--blown up then there is a vertex $\s$ of $S$ such that
$\conjugacyseq(\s)$ can be $\groupo$--blown up in $\group_\s$.\qed
\end{lemma}

\begin{lemma}\label{l:I-simplifying}
Let $\group$ be a group and let $S$ be a $\group$--tree. If $S$ can be
$\groupo$--I--simplified then $S$ has a valence one vertex with
stabilizer $\groupo$ whose incident edge has isomorphic
stabilizer.\qed
\end{lemma}

\begin{remark}
Recall
that a II--simplification has no effect on $S$.
\end{remark}

\section{Algorithmic Results} \label{s:algorithms}
\subsection{More labeled graphs}
A map of labeled graphs $\map\co \graph_1\to\graph_2$ is a {\it morphism} if 
\begin{itemize}
\item
the induced map between universal covers is a morphism;
and 
\item
$\map$ is label-preserving, ie, the following diagram commutes.
$$
\begindc{\commdiag}[40]
\obj(2,1)[rose]{$R_\B$}
\obj(1,2)[graph1]{$\graph_1$}
\obj(3,2)[graph2]{$\graph_2$}
\mor{graph1}{graph2}{$\map$}
\mor{graph1}{rose}{$\proj(\graph_1)$}[\atright,\solidarrow]
\mor{graph2}{rose}{$\proj(\graph_2)$}
\enddc
$$
\end{itemize}
It is {\it\proper}\ 
if this induced map is also \proper.
Stallings introduced
labeled graphs into the study of free groups. The next lemma is key.

\begin{lemma}[Stallings \cite{st:folding}]\label{l:stallings}
An immersion of labeled graphs induces an injection of
fundamental groups.\qed
\end{lemma}

An {\it edge path} in a labeled graph $\graph$ is a strict morphism
$I\to\graph$ where $I$ is an oriented compact interval.  If $I$ is a
point then the edge path is {\it trivial}. A nontrivial edge path may
be identified with a sequence of oriented edges
$\edge_0^{}\cdots\edge_m^{}$ where, for $1\le i\le m$,
$\partial_1\edge_{i-1}^{}=\partial_0 \edge_i^{}$.
The product of edge paths $\edgepath_1$ and $\edgepath_2$ is denoted
$\edgepath_1\edgepath_2$.  An edge path is {\it closed} if the initial
and terminal vertices of $I$ have the same image.

A {\it loop} in $\graph$ is a \proper\ morphism $S\to\graph$ where $S$
is an oriented circle.  A loop may be represented by a cyclic sequence
of edges of $\graph$. An oriented edge of $\graph$ is {\it crossed} by
a path or a loop if appears in the edge sequence representing the
path or loop.  The graph $\graph$ is {\it tight} if its labeling
function $\proj(\graph)\co \graph\to R_\B$ is an immersion. We record a
simple corollary of Lemma~\ref{l:stallings}.

\begin{corollary}\label{c:stallings}
If $\map\co I\to \graph$ is a tight non-trivial edge
path then the element of $\pi_1(R_\B)$ represented by
$\proj(\graph)\circ\map$ is non-trivial.  Equivalently, if
$\edge_0\cdots\edge_m$ represents a tight non-trivial edge path and if
the label of $\edge_i$ is $\lab_i^{\delta_i}$ ($\delta_i=\pm 1$) then
$\lab_0^{\delta_0}\cdots\lab_m^{\delta_m}$ is non-trivial in $\f(\B)$.\qed
\end{corollary}

\subsection{Operations on graphs}
\subsubsection{Coring}\label{s:coring}
A {\it core graph} is a graph such that every edge is crossed by an
immersed loop. By Zorn's lemma, every graph $\graph$ has a unique
maximal core subgraph, its {\it core}, denoted $\core(\graph)$. A core
graph contains no valence 0 or 1 vertices. If $\graph$ is connected
and has finite fundamental group, then $\core(\graph)$ is finite. The
core of a tree is empty. If $\graph$ is labeled, then so is
$\core(\graph)$. In fact, $\core$ is a functor from the category of
labeled graphs and immersions to the category of labeled core graphs
and immersions. The map $\core(\graph)\subset\graph$ is natural with
respect to this functor. The conjugacy class $[[\subgroup]]$ of a
subgroup of $\f(\B)$ is uniquely represented by the core
$\graph(\subgroup)$ of the cover $R_{\B,H}$ of $R_\B$ corresponding
to $\subgroup$. The simple proof of the next lemma is left to the reader.

\begin{lemma}\label{l:loop}
Let $\edge$ be an edge of the labeled graph $\graph$.
\begin{itemize}
\item
Suppose that $\edge$ does not separate $\graph$. Then,
$\edge\subset\core(\graph)$ if and only if there is an immersed loop
crossing $\edge$ exactly once.
\item
Suppose that $\edge$ separates $\graph$. Then,
$\edge\subset\core(\graph)$ if and only if there is an immersed loop
crossing each of $\edge$ and $\edge^{-1}$ exactly
once.\qed
\end{itemize}
\end{lemma}

\subsubsection{Folds and tightening}
A morphism $\map\co \graph_1\to\graph_2$ of graphs is a {\it fold} if
the induced map between universal covers is a fold.  A fold induces a
surjection on the level of fundamental groups. It is a homotopy
equivalence unless the edges that are identified share both initial
and terminal vertices \cite{fh:coherence}.

A finite graph $\graph$ may be iteratively folded until it is
tight. If $\graph$ is not finite, then the direct limit of the system
of finite sequences of folds is well-defined. The result is the {\it
tightening} of $\graph$ and is denoted $\tight(\graph)$. Fix a base
vertex for $\graph$ (if $\graph$ is non-empty) and let $\subgroup$
denote the image
$(\proj(\graph))_\#(\pi_1(\graph))\subset\pi_1(R_\B)$.  Then,
$\proj(\graph)$ lifts to $\graph\to R_{\B,H}$. The graph
$\tight(\graph)$ may be identified with the image of this lift. In
fact, $\tight$ is a functor from the category of labeled graphs and
\proper\ morphisms to the category of tight labeled graphs and
immersions. The quotient map $\graph\to\tight(\graph)$ is natural with
respect to this functor.  More generally, if $\lab\in\B$ and if
$\graph$ is a labeled graph, then we define $\tight_\lab(\graph)$ as
above except that only edges labeled $\lab$ or $\lab^{-1}$ are folded.

\subsubsection{Applying an automorphism}
If $\graph$ is a labeled graph and $\auto\in\Aut(\f(\B))$, then
$\auto\graph$ is the labeled graph obtained by replacing each labeled
oriented edge $\edge$ of $\graph$ by the sequence of labeled oriented
edges $\auto\edge$. More precisely, if the oriented edge $\edge$ has
the label $\lab$ and if $\auto(\lab)=\word$ where $\word$ is a reduced
word of length $k$ in $\B$, then $\auto\edge$ is obtained from
$\edge$ by subdividing $\edge$ into $k$ subedges. The $i^{th}$ letter
of $\word$ has the form $\llab^\delta$ where $\llab\in\B$, and
$\delta=\pm 1$. The $i^{th}$ subedge of $\edge$ is given the label
$\llab$ and an orientation agreeing with that of $\edge$ if $\delta$
is positive and the opposite orientation otherwise.
The operation of applying the automorphism $\auto$ is a functor from
the category of labeled graphs and morphisms to itself. The
construction gives a cellular map $\auto\co \graph\to\auto\graph$ that is well
defined up to a homotopy rel vertices and that is natural (but not a
morphism).

\begin{lemma}\label{l:pullback}
If $\edgepath$ is an immersed edge path in $\auto\graph$, then there are
an immersed edge path $\hat\edgepath$ in $\graph$ represented by
$\edge_0\cdots\edge_m$, an initial edge subpath $\edgepath_0$ of
$\auto(\edge_0)$, and a terminal edge subpath $\edgepath_m$ of
$\auto(\edge_m)$ such that
\begin{enumerate}
\item
$\edgepath_0\not=\auto(\edge_0)$;
\item
$\edgepath_m\not=\auto(\edge_m)$; and
\item
$\auto(\hat\edgepath)$ is the immersed edge path
$\edgepath_0\edgepath\edgepath_m$.
\end{enumerate}
\end{lemma}

\begin{proof}
We may view $\auto\graph$ as being obtained from $\graph$ by
subdividing and relabeling. With this in mind, any immersed edge path
$\edgepath$ in $\auto\graph$ gives an immersed path $\edgepath'$ in
$\graph$ that may not have endpoints vertices. This path extends
uniquely to an immersed edge path $\hat\edgepath$ that is minimal with
respect to containing $\edgepath'$.
\end{proof}

\subsubsection{Collapsing edges}
If $\edge$ is an edge of the labeled graph $\graph$, then
$\map\co \graph\to\graph'$ is a collapse of $\edge$ if the induced map
between universal covers is the morphism collapsing a lift of
$\edge$. In this case, we denote $\graph'$ by
$\collapse_\edge(\graph)$. More generally, if $\edgeset$ is a set of
edges in $\graph$ then we may collapse each edge in $\edgeset$ to a
point and obtain $\collapse_\edgeset(\graph)$.  If
$\morphism\co \graph\to\graph'$ is a morphism, and if $\edgeset'$ is a
set of edges in $\graph'$, then there is an induced morphism
$\collapse_{\edgeset'}(\morphism)\co \collapse_{\morphism^{-1}(\edgeset')}(\graph)\to\collapse_{\edgeset'}(\graph')$.
To each edge $\edge'$ in $\collapse_\edgeset(\graph)$, we may
associate the unique edge $\edge$ of $\graph$ such that
$\collapse_\edgeset(\edge)=\edge'$. The proof of the next lemma
is left to the reader.

\begin{lemma}[See \cite{fh:coherence}]\label{l:collapsing}
\begin{enumerate}
\item\label{i:collapsing 1} The quotient map
$\graph\to\collapse_\edgeset(\graph)$ induces a surjection of
fundamental groups.
\item \label{i:collapsing 2} If
$\edge'\subset\core(\collapse_\edgeset(\graph))$ then
$\edge\subset\core(\graph)$.
\end{enumerate}
\end{lemma}

\begin{remark}\label{r:collapsing and tightening commute}
If $\edgeset$ is the set of edges labeled $\lab$, then the operations
$\tight_\lab$ and $\collapse_\edgeset$ commute.
\end{remark}

\subsubsection{Operations on tight labeled core graphs}
If $\graph$ is a tight labeled core graph and if
$\auto\in\Aut(\f(\B))$, then $\auto_\#\graph$ is the labeled core
graph $\core(\tight(\auto\graph))$ obtained by coring the tightening
of $\auto\graph$. If $\graph$ represents $[[\subgroup]]$, then
$\auto_\#\graph$ represents $[[\auto\subgroup]]$.

\subsubsection{Sequences} 
All the above notions extend to sequences of labeled graphs. For
example, if $\graphseq=\{\graph_k\}$ is a sequence of labeled
graphs, then a path in $\graphseq$ is a path $I\to G_{k_0}$ for some
choice $k_0$ of $k$, $\auto_\#\graphseq$ denotes
$\{\auto_\#\graph_k\}$, etc.

A sequence $\conjugacyseq$ of conjugacy classes of subgroups of $\f(\B)$
is uniquely represented by $\graph(\conjugacyseq)$.  In the following
definitions, $\graphseq$ is a sequence of labeled graphs. For a
labeled graph $\graph$, $\abs(\lab,\graph)$ is the number of oriented
edges of $\graph$ with label $\lab$. If $\graphseq=\{\graph_i\}$ is a
sequence of labeled graphs then
$\abs(\lab,\graphseq)$ is the sum of the $\abs(\lab,\graph_i).$

\subsection{Elementary Whitehead automorphisms}
A reference for this section is \cite{sk:gersten}. An {\it extended
permutation of $\f(\B)$} is an automorphism of $\f(\B)$ induced by a permutation
of $\B^{\pm 1}$. An {\it elementary Whitehead automorphism} is an
automorphism $\auto$ of $\f(\B)$ that is either an extended
permutation or has the following form. There is an element
$\lab\in\B^{\pm 1}$ and a subset $A$ of $\B^{\pm
1}\setminus\{\lab^{\pm 1}\}$ such that
\begin{itemize}
\item
if $\llab\in A\setminus A^{-1}$ then $\auto(\llab)=\lab\llab$;
\item
if $\llab\in A\cap A^{-1}$ then $\auto(\llab)=\lab\llab \lab^{-1}$; and
\item
if $\llab\not\in A\cup A^{-1}$ then $\auto(\llab)=\llab$.
\end{itemize}
We call $\lab$ {\it the distinguished label of $\auto$.}

\begin{remark}\label{r:correspondence}
Let $\graphseq$ be a sequence of labeled graphs and let $\auto$ be an
elementary Whitehead automorphism with distinguished label
$\lab$. There is a 1--1 correspondence between the set of edges of
$\graphseq$ not labeled $\lab$ and the edges of $\auto\graphseq$ not
labeled $\lab$. In $\auto\graphseq$, there are {\it old} and {\it
new} edges labeled $\lab$. The terminal vertex of each new edge has
valence 2 and the other incident edge is not labeled $\lab$. Such a
valence 2 vertex is {\it new}\,; other vertices are {\it old}.
The subgraph of $\auto\graphseq$ consisting of new edges is a forest
each component of which is a cone over a set of new vertices with base
an old vertex. All edges of the cone have initial vertex the base. See
\figref{f:auto}.
\begin{figure}[ht!]\anchor{f:auto}
\cl{\scalebox{0.75}{\includegraphics{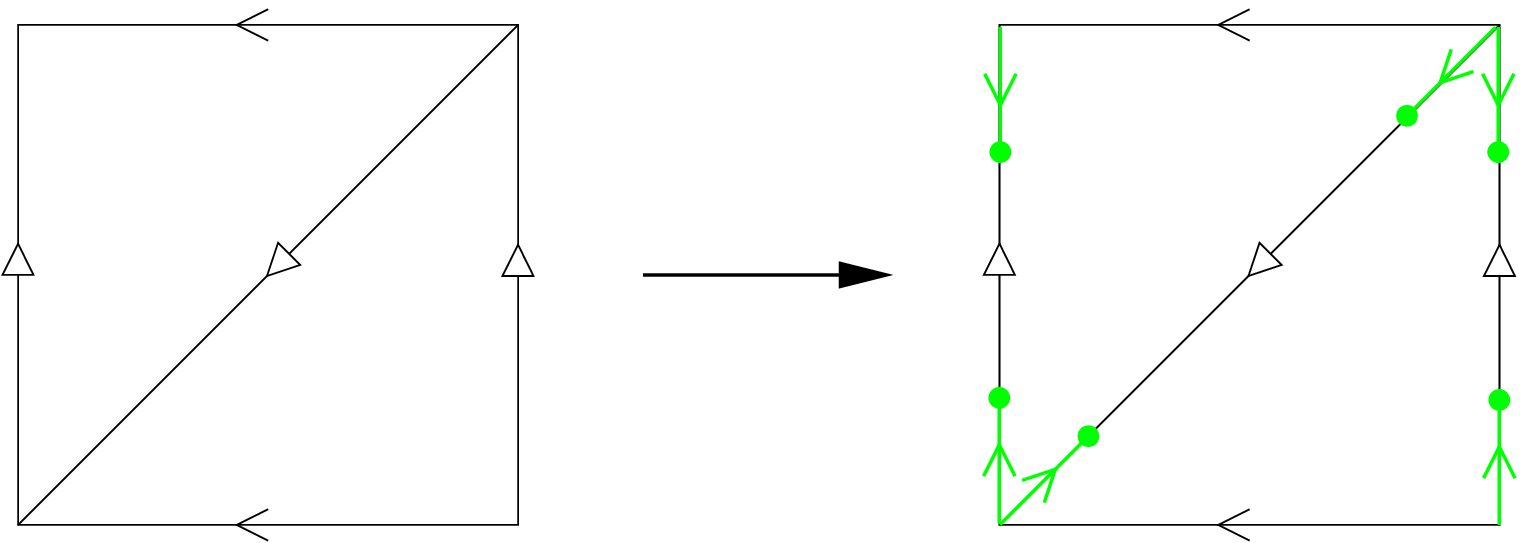}}}
\caption{$\auto(a)=a$, $\auto(b)=aba^{-1}$}\label{f:auto}
\end{figure}
\end{remark}

\begin{remark}\label{r:whitehead}
For $\auto\in\Aut(\f(\B))$, the sequence of folds needed to tighten
$\auto R_\B$ algorithmically gives a factorization of $\auto$ as a
product of elementary Whitehead automorphisms.
\end{remark}

The next lemma is a consequence of Step~1 of the proof of the proposition on
page~455 of \cite{bf:bounding}.
\begin{lemma}\label{l:fold}
Let $\map\co \graph_0\to\graph_1$ be a \proper\ morphism of labeled
graphs that is surjective on the level of fundamental groups. Then,
there is a fold $\map'\co \graph_0\to\graph'$ such that $\map$ factors as
$$\graph_0\overset{\map'}{\to}\graph'\to\graph_1.\qed$$
\end{lemma}

\begin{lemma}\label{l:path}
Let $\map\co \graph_0\to\graph_1$ be a \proper\ morphism of labeled
graphs that is surjective on the level of fundamental groups. Suppose
that, for some $\lab\in\B$,
$\abs(\lab,{\graph_1})<\abs(\lab,{\graph_0})$. Then, there are
\proper\ morphisms making the following diagram commute
$$
\begindc{\commdiag}[40]
\obj(1,2)[I]{$I$}
\obj(2,2)[graph0]{$\graph_0$}
\obj(1,1)[T]{$T$}
\obj(2,1)[graph1]{$\graph_1$}
\mor{I}{graph0}{$\edgepath$}
\mor{I}{T}{$\map'$}[\atright,\solidarrow]
\mor{T}{graph1}{}
\mor{graph0}{graph1}{$\map$}
\enddc
$$
where
\begin{itemize}
\item
where $\edgepath$ is an immersed edge path represented by
$\edge_0\edge_1\cdots\edge_m$;
\item
$T$ is a labeled tree;
\item
$\edge_0$ and $\edge_m^{-1}$ are labeled by $\lab'=\lab$ or $\lab^{-1}$;
\item
for $0<i<m$, $\edge_i$ is not labeled by $\lab$ or $\lab^{-1}$; and 
\item
$\map'(\edge_0)=\map'(\edge_m^{-1})$.
\end{itemize}
\end{lemma}

\begin{proof}
Since $\abs(\lab,{\graph_1})<\abs(\lab,{\graph_0})$ there are distinct
edges $\edge$ and $\edge'$ in $\graph_0$ each labeled $\lab$ that are
identified under $\map$. Consider the lift
$\tilde\map\co \tilde\graph_0\to\tilde\graph_1$ to universal
covers. Because $\map$ induces a surjection on the level of
fundamental groups, there are lifts $\tilde\edge$ and $\tilde\edge'$
of $\edge$ and $\edge'$ to $\tilde\graph_0$ that are identified under
$\tilde\map$. Choose $\tilde\edge$ and $\tilde\edge'$ with this
property so that the subtree $I$ they span has minimal diameter (with
respect to the edge metric). The edge path $\edgepath\co I\to\graph_0$
is the restriction of the first covering projection. The edge path
$\edgepath$ factors as $I\to T=\tilde\map(I)\to\graph_1$ where the
first factor is induced by the restriction of $\tilde\map$ to $I$ and
the second factor is the restriction of the second covering
projection.
\end{proof}

\begin{lemma}\label{l:collapse iso}
Let $\graphseq$ be a sequence of tight labeled graphs and let
$\auto\in\Aut(\f(\B))$ be an elementary Whitehead automorphism with
distinguished label $\lab$. Let $\edgeset$ (respectively $\edgeset'$) be
the set of edges of $\graph$ (respectively $\tight(\auto\graph)$) that
are labeled $\lab$. Then, the following diagram commutes and the lower
horizontal arrow is an isomorphism.
$$
\begindc{\commdiag}[40]
\obj(1,2)[graph]{$\graphseq$}
\obj(4,2)[tightautograph]{$\tight(\auto\graphseq)$}
\obj(1,1)[collapsegraph]{$\collapse_\edgeset(\graphseq)$}
\obj(4,1)[collapsetightautograph]{$\collapse_{\edgeset'}(\tight(\auto\graphseq))$}
\mor{graph}{tightautograph}{}
\mor{graph}{collapsegraph}{}
\mor{tightautograph}{collapsetightautograph}{}
\mor{collapsegraph}{collapsetightautograph}{}
\enddc
$$ In particular, there is a natural 1--1 correspondence between
edges of $\graphseq$ not labeled $\lab^{\pm 1}$ and edges of
$\tight(\auto\graphseq)$ not labeled $\lab^{\pm 1}$.
\end{lemma}

\begin{proof}
Let $\edgeset''$ be the set of edges of $\auto\graph$ that are labeled
$\lab$. We have a commuting diagram
$$
\begindc{\commdiag}[40]
\obj(1,2)[graph]{$\graphseq$}
\obj(4,2)[autograph]{$\auto\graphseq$}
\obj(8,2)[tightautograph]{$\tight(\auto\graphseq)$}
\obj(1,1)[collapsegraph]{$\collapse_\edgeset(\graphseq)$}
\obj(4,1)[collapseautograph]{$\collapse_{\edgeset''}(\auto\graphseq)$}
\obj(8,1)[collapsetightautograph]{$\collapse_{\edgeset'}(\tight(\auto\graphseq))$}
\mor{graph}{autograph}{}
\mor{autograph}{tightautograph}{}
\mor{graph}{collapsegraph}{}
\mor{autograph}{collapseautograph}{}
\mor{tightautograph}{collapsetightautograph}{}
\mor{collapsegraph}{collapseautograph}{}
\mor{collapseautograph}{collapsetightautograph}{}
\enddc
$$ It is clear that the lower left horizontal arrow is an isomorphism
and that the lower right arrow is \proper\ and surjective. 

In order to obtain a contradiction, assume
$$\collapse_{\edgeset''}(\auto\graphseq)\to\collapse_{\edgeset'}(\tight(\auto\graphseq))$$
is not injective. Since this map is $\pi_1$--surjective, by
Lemma~\ref{l:fold} there are two edges not labeled $b$ or $b^{-1}$
with the same image. It follows that there are two edges not labeled
$b$ or $b^{-1}$ with the same image under
$\auto\graphseq\to\tight(\auto\graphseq).$ Using Lemma~\ref{l:path}
and taking a subpath if necessary, there is an immersed edge path
$\edgepath\co I\to\auto\graphseq$ represented by
$\edge_0\edge_1\cdots\edge_m$ such that
\begin{itemize}
\item
the label of $\edge_0$ and $\edge_m^{-1}$ is $\lab'\not=\lab^{\pm 1}$;
\item
the label of $\edge_i$ is $\lab^{\pm 1}$ for all $0<i<m$; and
\item
$I\to\auto\graphseq\to\tight(\auto\graphseq)$ factors through a tree.
\end{itemize}

If $\hat\edgepath$ is the immersed edge path in $\graphseq$ determined by
$\edgepath$ as in Lemma~\ref{l:pullback} then
\begin{itemize}
\item
the label of $\hat\edge_0$ and $\hat\edge_{\hat m}^{-1}$ is $\lab'$;
and
\item
$\hat\edge_i$ is labeled $\lab^{\pm 1}$ for $1<i<\hat m$.
\end{itemize}
Since $\hat\edgepath$ is an immersion, all the $\hat\edge_i$, $1<i<\hat
m$, are consistently oriented. It is easy to see then that
$I\to\auto\graphseq\to\tight(\auto\graphseq)$ cannot factor
through a tree, contradiction.
\end{proof}

\begin{lemma}\label{l:core labels}
Let $\graphseq$ be a sequence of tight labeled core graphs and let
$\auto$ be an elementary Whitehead automorphism with distinguished
label $\lab$, let $\edge$ be an edge of $\graphseq$ not labeled
$\lab^{\pm 1}$, and let $\edge'$ be the corresponding edge in
$\tight(\auto\graphseq)$. Then, $\edge'$ is in
$\auto_\#\graphseq=\core(\tight(\auto\graphseq))$. In particular,
there is a natural 1--1 correspondence between edges of\, $\graphseq$ not
labeled $\lab^{\pm 1}$ and edges of $\auto_\#\graphseq$ not labeled
$\lab^{\pm 1}$.
\end{lemma}

\begin{proof}
Suppose that $\edge$ separates (respectively does not separate) its
component. By Lemma~\ref{l:loop}, there is an immersed loop
$\map\co S\to\graphseq$ crossing $\edge$ (respectively crossing $\edge'$
and $\edge'^{-1}$ each) exactly once. It follows from
Lemma~\ref{l:collapse iso} that the immersed loop $\tight(\auto\map)$
crosses $\edge'$ (respectively $\edge'$ and $\edge '^{-1}$ each)
exactly once. By Lemma~\ref{l:loop}, $\edge'$ is contained in
$\auto_\#\graphseq$.
\end{proof}

\subsection{Complexity}
If $\conjugacyseq$ is a finite sequence of conjugacy
classes of finitely generated subgroups of $\f(\B)$ and if $\lab\in\B$,
then $\abs(\lab,\conjugacyseq)$ is the number of edges in
$\graph(\conjugacyseq)$ that are labeled with
$\lab$. The {\it complexity of $\conjugacyseq$}, denoted
$\complexity(\conjugacyseq)$, is the number of edges in
$\graph(\conjugacyseq)$ or equivalently $\sum_{\lab\in\B}
\abs(\lab,\conjugacyseq)$.

We will also need a finer measure of complexity of
$\conjugacyseq$. Define the {\it lexity} of $\conjugacyseq$, denoted
$\lex(\conjugacyseq)$, to be the sequence of non-negative integers
$\{\abs(\lab,{\graph(\conjugacyseq)})\}_{\lab\in\B}$ arranged in
non-decreasing order. The set $\lexspace$ of non-decreasing sequences
of non-negative integers is well-ordered lexicographically. Let
$\minlex(\conjugacyseq)$ denote
$\min_{\lab\in\B}\{\abs(\lab,{\graph(\conjugacyseq)})\}$.

\begin{lemma}\label{l:c vs abs}
Let $\conjugacyseq$ be a finite sequence of conjugacy classes of
finitely generated subgroups of $\f(\B)$ and let $\auto$ be an elementary
Whitehead automorphism. Then
$\complexity(\auto\conjugacyseq)<\complexity(\conjugacyseq)$ if and
only if $\lex(\auto\conjugacyseq)<\lex(\conjugacyseq)$. If further
$\auto$ has distinguished label $\lab$, then
$\lex(\auto\conjugacyseq)\le\lex(\conjugacyseq)$ if and only if
$\abs(\lab,\auto\conjugacyseq)\le\abs(\lab,\conjugacyseq)$ with
equality if and only if
$\abs(\lab,\auto\conjugacyseq)=\abs(\lab,\conjugacyseq)$.
\end{lemma}

\begin{proof}
Since extended permutations preserve both $\complexity$ and $\lex$, we
may suppose that $\auto$ has distinguished label $\lab\in\B$. It
follows from Lemma~\ref{l:core labels} that if
$\lab'\not=\lab^{\pm 1}$, then the number of times that $\lab'$ appears in
$\graph(\auto\conjugacyseq)$ is the same as the number of times that
$\lab'$ appears in $\graph(\conjugacyseq)$.
\end{proof}

\subsection{Gersten's Theorem}\label{s:gersten}
Let $\conjugacyseq$ be a finite sequence of conjugacy classes of
finitely generated subgroups of $\f(\B)$. If
$\complexity(\conjugacyseq)=\min\{\complexity(\auto\conjugacyseq)\mid\auto\in\Aut(\f(\B))\}$
then $\conjugacyseq$ is a {\it Gersten representative for the orbit
$\Aut(\f(\B))\conjugacyseq$}. We also write that $\conjugacyseq$ is a Gersten
representative for any element of the orbit. Since $\conjugacyseq$ is
an element of this orbit, we often simply write that $\conjugacyseq$
is a Gersten representative. A finite set of generators for a
representative $\subgroup\in\conjugacy$ for each
$\conjugacy\in\conjugacyseq$ is a {\it finite generating system} for
$\conjugacyseq$.  S\,M~Gersten \cite{sg:whitehead}\cite{sk:gersten} gave
an algorithm that when input a finite generating system for $\conjugacyseq$
outputs the finite set of Gersten representatives for $\conjugacyseq$.

\begin{theorem}{\rm\cite{sg:whitehead},\cite{sk:gersten}}\label{t:gersten}
\begin{enumerate}
\item\label{i:gersten 1} If $\conjugacyseq$ is a finite sequence of
conjugacy classes of finitely generated subgroups of $\f(\B)$ that is not
a Gersten representative, then there is an elementary
Whitehead automorphism $\auto$ such that
$\complexity(\auto\conjugacyseq)<\complexity(\conjugacyseq)$.
\item\label{i:gersten 2}
If $\conjugacyseq$ and $\conjugacyseq'$ are Gersten representatives for
$\conjugacyseq$, then there is a finite sequence $\{\auto_k\}_{k=1}^m$ of
elementary Whitehead automorphisms and a sequence
$$\{\conjugacyseq=\conjugacyseq_0,\conjugacyseq_1,\cdots,\conjugacyseq_m=\conjugacyseq'\}$$
of Gersten representatives such that
$\auto_{k}\conjugacyseq_{k-1}=\conjugacyseq_{k}$ for $1\le k\le m$.
\end{enumerate}
\end{theorem}

\begin{corollary}\label{c:gersten 1}
If $\conjugacyseq$ is a finite sequence of conjugacy classes of finitely
generated subgroups of $\f(\B)$, then there is an algorithm that when input a
finite generating system for $\conjugacyseq$ outputs a Gersten
representative for $\conjugacyseq$.
\end{corollary}

\begin{proof}
Iteratively apply elementary Whitehead automorphisms to $\conjugacyseq$
until complexity cannot be decreased. The resulting sequence is a
Gersten representative by Theorem~\ref{t:gersten}(\ref{i:gersten 1}).
\end{proof}

\begin{corollary}
Let $\conjugacyseq$ be a finite sequence of conjugacy classes of
finitely generated subgroups of $\f(\B)$. Then, there is an algorithm that
when input a finite generating system for $\conjugacyseq$ outputs the
finitely many Gersten representatives of $\conjugacyseq$.
\end{corollary}

\begin{proof}
By Corollary~\ref{c:gersten 1}, we may assume that $\conjugacyseq$ is a
Gersten representative. Consider the graph whose vertices are finite
sequences of conjugacy classes of finitely generated subgroups of $\f(\B)$
of complexity equal to $\complexity(\conjugacyseq)$, and where two
vertices $\conjugacyseq_1$ and $\conjugacyseq_2$ are connected by an
edge if there is an elementary Whitehead automorphism $\auto$ such
that $\auto\conjugacyseq_1=\conjugacyseq_2$. By
Theorem~\ref{t:gersten}(\ref{i:gersten 2}), the component of
this graph containing $\conjugacyseq$ has vertices that are precisely
the Gersten representatives of $\conjugacyseq$.
\end{proof}

\subsection{Consequences of Lemma~\ref{l:c vs abs}}
In this section we show that simplifications can be detected using
Gersten representatives. Throughout this section, $\conjugacyseq$ is a
finite sequence of conjugacy classes of finitely generated subgroups
of $\f(\basis)$. 

\begin{lemma} \label{l:min lex 0}
The following are equivalent.
\begin{enumerate}
\item \label{i:min lex 0 1}
There is an $\auto\in\Aut(\f(\B))$ such that $\minlex(\auto\conjugacyseq)=0$.
\item\label{i:min lex 0 2} 
For some (any) Gersten representative
$\conjugacyseq'$ of $\conjugacyseq$, $\minlex(\conjugacyseq')=0$.
\end{enumerate}
\end{lemma}

\begin{proof}
By Lemma~\ref{l:c vs abs}, for any $\auto\in\Aut(\f(\B))$ and any Gersten
representative $\conjugacyseq'$ of $\conjugacyseq$,
$\lex(\auto\conjugacyseq)\ge\lex(\conjugacyseq')$. Thus, (\ref{i:min
lex 0 1})$\iff$(\ref{i:min lex 0 2}). 
\end{proof}

\begin{lemma}\label{l:fewer labels}
Let $\B'\subset\B$ and suppose that $\conjugacyseq'$ is a
finite sequence of conjugacy classes of finitely generated subgroups
of $\f(\B')$. The following are equivalent.
\begin{enumerate}
\item
$\conjugacyseq'$ is a Gersten representative with respect to $\B'$.
\item
$\conjugacyseq'$ is a Gersten representative with respect to $\B$.
\end{enumerate}
\end{lemma}

\begin{proof}
That (2)$\implies$(1) is clear. Suppose (1), but not (2). By
Lemma~\ref{l:c vs abs} there is then an elementary Whitehead automorphism
$\auto\in\Aut(\f(\B))$ with distinguished label
$\lab\notin(\B')^{\pm 1}$ such that
$\abs(\lab,\auto\conjugacyseq')<\abs(\lab,\conjugacyseq')=0$,
contradiction.
\end{proof}

Recall that the terms {\it visibly blown up}, {\it visibly unpulled},
{\it visibly unkilled}, {\it visibly uncleaved}, and {\it visibly
simplified} were defined in Definition~\ref{d:visible}.

\begin{lemma} \label{l:blow up}
The following are equivalent.
\begin{enumerate}
\item \label{i:partition 1} There is an $\auto\in\Aut(\f(\B))$ such that
$\graph(\auto\conjugacyseq)$ can be visibly blown up.
\item \label{i:partition 2} For some (every) Gersten representative
$\conjugacyseq'$ of $\conjugacyseq$, $\graph(\conjugacyseq')$ can be
visibly blown up.
\end{enumerate}
\end{lemma}

\begin{proof}
The lemma will follow from:
\medskip

{\bf Claim}\qua If $\graph(\conjugacyseq)$ can be visibly blown
up and if $\auto$ is an elementary Whitehead automorphism with
$\lex(\auto\conjugacyseq)\le\lex(\conjugacyseq)$ then either
$\graph(\alpha\conjugacyseq)$ can be visibly blown up or there is an
automorphism $\auto'$ with
$\lex(\auto'\conjugacyseq)<\lex(\conjugacyseq)$ such that
$\graph(\auto'\conjugacyseq)$ can be visibly blown up.
\medskip

We now prove the claim. Suppose $\B=\B'\sqcup\B''$ is a
non-trivial partition such that each element of $\conjugacyseq$ has a
representative in either $\f(\B')$ or $\f(\B'')$.  If $\auto$ is
an extended permutation, then the claim is clear.

Alternatively, let $\lab\in\B^{\pm 1}$ be the distinguished label
of $\auto$ and suppose without loss that $\lab\in(\B')^{\pm
1}$. Let $\auto'$ be the automorphism that agrees with $\auto$ on
$\B'$ and that is the identity on $\B''$. It is clear that
$\graph(\auto'\conjugacyseq)$ can be visibly blown up. Now,
$\graph(\auto\conjugacyseq)$ can be visibly blown up if and only if, for each
$\conjugacy\in\conjugacyseq$ with a representative in $\f(\B'')$,
$\alpha\conjugacy$ has a representative in $\f(\B'')$. By
Lemma~\ref{l:c vs abs}, this occurs if and only if
$\abs(\lab,\alpha\conjugacy)>0$ for such $\conjugacy$ and this occurs
if and only if
$\abs(\lab,\auto'\conjugacyseq)<\abs(\lab,\auto\conjugacyseq)$.
\end{proof}

\begin{lemma} \label{l:min lex 1}
Suppose that for some (any) Gersten representative $\conjugacyseq'$ of
$\conjugacyseq$ we have $\minlex(\conjugacyseq')\not= 0$. Then, the
following are equivalent.
\begin{enumerate}
\item \label{i:min lex 1 1}
There is an $\auto\in\Aut(\f(\B))$ such that $\minlex(\auto\conjugacyseq)=1$.
\item \label{i:min lex 1 2}
For some (any) Gersten representative $\conjugacyseq'$ of $\conjugacyseq$,
$\minlex(\conjugacyseq')=1$.
\item \label{i:min lex 1 3}
$\conjugacyseq$ can be either visibly unkilled or visibly unpulled.
\end{enumerate}
\end{lemma}

\begin{proof}
In the presence of $\minlex(\conjugacyseq')\not= 0$, (\ref{i:min lex 1
1})$\iff$(\ref{i:min lex 1 2}) by Lemma~\ref{l:c vs abs}. Suppose that
(\ref{i:min lex 1 2}) holds. There are two cases: there is a label
$\lab\in\B$ that appears exactly once in $\graph(\conjugacyseq')$ and
(a) the edge labeled $\lab$ separates its component and (b) edge
labeled $\lab$ does not separate its component.  It is an easy
exercise to show that in case~(a) $\graph(\conjugacyseq)$ can be
visibly unkilled and in case~(b) $\graph(\conjugacyseq)$ can be
visibly unpulled.

Suppose that (\ref{i:min lex 1 3}) holds. Then there is free
factorization $\f(\B)=F^1*_{\langle 1\rangle}=F^1*\langle t\rangle$ as in
Definition~\ref{d:algebraic unkilling} or Definition~\ref{d:algebraic
unpulling}. Let $\B^1$ be a basis for $F^1$. Choose an
$\auto\in\Aut(\f(\B))$ so that $\auto(\B^1\sqcup\{t\})=\B$. Then
$\minlex(\auto\conjugacyseq)=1$.
\end{proof}

The next lemma will be used to prove Lemma~\ref{l:split}.

\begin{lemma}\label{l:sublemma}
Suppose that
\begin{enumerate}
\item
$\minlex(\conjugacyseq)>1$;
\item
$\B=\B'\sqcup\B''$ is a non-trivial
partition;
\item
$\graph(\conjugacyseq)=\graph(\conjugacyseq')\vee\graph(\conjugacyseq'')$
where $\conjugacyseq'$ is a Gersten representative in $\f(\B')$
and $\conjugacyseq''$ is a Gersten representative in $\f(\B'')$; and
\item
$\auto\in\Aut(\f(\B))$ is an elementary Whitehead automorphism with
distinguished label $\lab\in(\B')^{\pm 1}$ such that
$\lex(\auto\conjugacyseq)\le\lex(\conjugacyseq)$.
\end{enumerate}
Then,
\begin{itemize}
\item
$\auto\conjugacyseq''=\conjugacyseq''$;
\item
$\auto\conjugacyseq'$ is a Gersten representative for
$\conjugacyseq'$; and
\item
$\graph(\auto\conjugacyseq)=\graph(\auto\conjugacyseq')\vee\graph(\conjugacyseq'')$.
\end{itemize}
In particular, $\conjugacyseq$ satisfying (1), (2), and (3) is a Gersten
representative.
\end{lemma}

\begin{proof}
Note:
\begin{itemize}
\item
$\graph(\auto\conjugacyseq)=\tight\big(\graph(\auto\conjugacyseq')\vee
I\vee\graph(\auto\conjugacyseq'')\big)$ for some labeled graph $I$
homeomorphic to a compact interval. This follows because in tightening
$\auto\graph(\conjugacyseq)$, we can tighten
$\auto\graph(\conjugacyseq')$ and $\auto\graph(\conjugacyseq'')$ first. 
\item
The subgraph of $\graph(\auto\conjugacyseq'')$ consisting of edges
labeled $\lab$ is a tree whose components are single (non-loop)
edges. This follows from Remark~\ref{r:correspondence} and the
following commutative diagram where $\edgeset$ is the set of edges of
$\graph(\auto\conjugacyseq'')$ that are labeled $\lab$.
$$
\begindc{\commdiag}[40]
\obj(1,2)[graph]{$\graph(\conjugacyseq'')$}
\obj(5,2)[autograph]{$\graph(\auto\conjugacyseq'')$}
\obj(3,1)[graphagain]{$\graph(\conjugacyseq'')$}
\mor{graph}{autograph}{}
\mor{graph}{graphagain}{=}[\atright,\solidarrow]
\mor{autograph}{graphagain}{$\collapse_\edgeset$}
\enddc
$$
\end{itemize}
It follows that in tightening $\graph(\auto\conjugacyseq')\vee
I\vee\graph(\auto\conjugacyseq'')$ at most one edge of
$\graph(\auto\conjugacyseq')$ folds with an edge of
$\graph(\auto\conjugacyseq'')$. Hence,
\begin{align*}
\complexity(\conjugacyseq')+\complexity(\conjugacyseq'')&=\complexity(\conjugacyseq)\ge
\complexity(\auto\conjugacyseq)\\ &\ge \complexity(\auto\conjugacyseq')+\complexity(\auto\conjugacyseq'')-1\ge\complexity(\conjugacyseq')+\complexity(\conjugacyseq'')-1.
\end{align*}
Thus, either
$\complexity(\auto\conjugacyseq'')=\complexity(\conjugacyseq'')$ or
$\complexity(\auto\conjugacyseq'')=\complexity(\conjugacyseq'')+1$. In
the former case, we are done. The latter case cannot occur. Indeed,
otherwise $\abs(\lab,\auto\conjugacyseq'')=1$. But then, by
Lemmas~\ref{l:fewer labels} and \ref{l:min lex 1} ,
$\minlex(\conjugacyseq'')\le 1$ and hence $\minlex(\conjugacyseq)\le
1$, contradiction.
\end{proof}

\begin{remark}
Without (1), Lemma~\ref{l:sublemma} is false. Consider
$\B=\{a,b\}\sqcup\{c\}$, $H=\langle ab^{-1},c\rangle$, and $\auto(a)=a$,
$\auto(b)=b$, $\auto(c)=bc$.
\end{remark}

\begin{lemma}\label{l:split}
Suppose that $\minlex(\conjugacyseq')>1$ for some (any) Gersten
representative $\conjugacyseq'$ of $\conjugacyseq$. Then, the
following are equivalent.
\begin{enumerate}
\item\label{i:split 2}
There is $\auto\in\Aut(\f(\B))$ such that $\graph(\auto\conjugacyseq)$ can
be visibly \splited.
\item\label{i:split 3}
For some (every) Gersten representative $\conjugacyseq'$ of
$\conjugacyseq$, $\graph(\conjugacyseq')$ can be visibly \splited.
\end{enumerate}
\end{lemma}

\begin{proof} 
(\ref{i:split 3})$\implies$(\ref{i:split 2}) is clear. We now show
(\ref{i:split 2})$\implies$(\ref{i:split 3}).  Suppose that
$\graph(\auto\conjugacyseq)=\graph(\conjugacyseq')\vee\graph(\conjugacyseq'')$
where all labels of $\graph(\conjugacyseq')$ are in $\B'^{\pm 1}$
and all labels of $\graph(\conjugacyseq'')$ are in $\B''^{\pm 1}$
for some non-trivial partition $\B=\B'\sqcup\B''$. Choose
$\auto'\in\Aut(\f(\B'))$ and $\auto''\in\Aut(\f(\B''))$ such
that $\auto'\conjugacyseq'$ and $\auto''\conjugacyseq''$ are Gersten
representatives. Let $\auto\in\Aut(\f(\B))$ agree with $\auto'$ on
$\f(\B')$ and $\auto''$ on $\f(\B'')$. Then,
$\graph(\auto\conjugacyseq)=\graph(\auto'\conjugacyseq')\vee I'\vee
I''\vee\graph(\auto''\conjugacyseq'')$ where $I'$ and $I''$ are labeled
graphs homeomorphic to compact intervals, all labels of $I'$ are in
$(\B')^{\pm 1}$, and all labels of $I''$ are in $(\B'')^{\pm
1}$. We will now show that
$\graph(\auto'\conjugacyseq')\vee\graph(\auto''\conjugacyseq'')$ is also
a representative for $\conjugacyseq$. Suppose that $I''$ is not
trivial and that the edge $\edge$ of $I''$ with initial vertex in
$\graph(\auto'\conjugacyseq')\vee I'$ is labeled $\lab$. Let
$\auto_\lab$ be the elementary Whitehead automorphism that is
conjugation by $\lab$ on $\B'$ and the identity on
$\B''$. Then, $\graph(\auto_b\auto\conjugacyseq)$ is obtained from
$\graph(\auto'\conjugacyseq')\vee I'\vee
I''\vee\graph(\auto''\conjugacyseq'')$ by collapsing $\edge$. We may
continue until $I''$ and symmetrically $I'$ are trivial. It follows
from Lemma~\ref{l:sublemma} that the result is a Gersten
representative and that all Gersten representatives have this form.
\end{proof}

\section{Proof of the Main Theorem}\label{s:final}
\begin{proposition}\label{p:algorithm}
The algorithm of Section~\ref{ss:algorithm} is in fact an algorithm.
\end{proposition}

\begin{proof}
To detect a reduction, it is only necessary to be able to decide
algorithmically if a homomorphism $\auto\co \f(\B)\to \f(\B')$ between
free groups is an isomorphism. According to Stallings
\cite{st:folding}, $\auto$ is injective if and only if the rank of the
Stallings representative of $\auto F(\B)$ with respect to $\B'$ is
$|\B|$. It is surjective if and only if this Stallings representative
is $R_{\B'}$. Thus Step~1 is algorithmic.

Step~2 only depends on being able to find a Gersten representative and
this is algorithmic by Corollary~\ref{c:gersten 1}.

After reducing and $\trivialgroup$--simplifying, complexity has been
reduced where complexity is the sequence of ranks of conjugacy classes
of edge stabilizers viewed as an element of $\lexspace$. Therefore
this process stops.
\end{proof}

We are finally in a position to prove Theorem~\ref{t:main 2}.

\begin{proof}[Proof of Theorem~\ref{t:main 2}]
By Proposition~\ref{p:algorithm}, we may assume that $\G$ is
reduced. Let $S$ be the corresponding $G$--tree. We may also assume
that no edge stabilizer is trivial of $S$. (Otherwise, $\group$ has an
obvious free decomposition and we may work with the factors instead of
$\group$.) If $\group$ is freely decomposable, then by
Lemma~\ref{l:simplify characterization}, there is a vertex $\bsv\in V$
such that $\conjugacyseq(\bsv)$ can be $\trivialgroup$--simplified. In
particular, by definition there is a basis $\B$ for $\group_\bsv$
with respect to which $\graph(\conjugacyseq(\bsv))$ may be visibly
simplified. Since $\Aut(\f(\B))$ acts transitively on bases,
$\graph_\gersten(\conjugacyseq(\bsv))$ may be visibly simplified by
Lemmas~\ref{l:blow up}, \ref{l:min lex 1}, and \ref{l:split}.
\end{proof}

We end with a few questions.

\begin{question}
Is there an algorithm to decide if the fundamental group of a finite
graph of finite rank free groups is a surface group?
\end{question}

\begin{question}
Is there an algorithm to decide if the fundamental group of a finite
graph of finite rank free groups splits over $\mathbb Z$?
\end{question}

\begin{question}
Is there an algorithm to find the $JSJ$--decomposition of the
fundamental group of a finite graph of finite rank free groups?
\end{question}

One can't hope to go too far in this direction since according to
C.~Miller \cite{cfm:decision} the isomorphism
problem for finite graphs of finite rank free groups is unsolvable.

\section*{Appendix}
\addcontentsline{toc}{section}{Appendix}

\appendix
\section{Bookkeeping} 
In this section, details are provided as to how we record a finite
graph of finite rank free groups $\G$ and how that data changes under
a simplification. See Section~\ref{s:graphs of groups} for notation.

A finite graph of finite rank free groups $\G$ is given by the
following data:
\begin{itemize}
\item
for each $\bse\in\hat E$, a basis $\B_{\bse}$ for $\group_{\bse}$ such
that $\B_{\bse}=\B_{\bse^{-1}}$;
\item 
for $\bsv\in V$, a basis $\B_{\bsv}$ for $\group_{\bsv}$; and
\item
$\vec{w}=\{\vec{w}_{\bse}\}_{\bse\in\hat E}$ where
$\vec{w}_{\bse}=\{w_{{\bse},b}\}_{b\in \B_{\bse}}$ is a sequence of
reduced words in $\B_{\bsv}^{\pm 1}$ representing
$\{\varphi_{\bse}(b)\}_{b\in\B_{\bse}}$.
\end{itemize}

The Stallings algorithm referred to in Section~\ref{s:stallings} can
be used to decide if a sequence $\vec{w}_{\bse}$ determines a
monomorphism.  Indeed, check if the rank of $\graph_{\stallings}$
obtained from $\vec{w}_{\bse}$ is equal to the rank of
$\group_{\bse}$. 

\begin{definition}\label{d:conjugate word sequence}
If $\{\vec{w}_{\bse}\}_{{\bse}\in \hat E}$ is a sequence as above,
then we say another sequence $\{\vec{w}^\output_{\bse}\}_{{\bse}\in\hat E}$ is
{\it conjugate} to $\{\vec{w}_{\bse}\}_{{\bse}\in \hat E}$, written
$\{\vec{w}^\output_{\bse}\}_{{\bse}\in \hat
E}\sim\{\vec{w}_{\bse}\}_{{\bse}\in\hat E}$, if there are
$\{\psi_{\bse}=\psi_{{\bse}^{-1}}\in\Aut(\group_{\bse})\}_{{\bse}\in\hat
E}$, $\{\psi_{\bsv}\in\Aut(\group_{\bsv})\}_{{\bsv}\in V}$, and
$\vec{h}=\{h_{\bse}\in \group_{\bsv}\}_{{\bse}\in\hat E }$
so that $\vec{w}^\output_{\bse}$ represents $\{\psi_{\bsv}
\circ i_{h_{\bse}} \circ\varphi_{\bse}
\circ \psi_{\bse}(b)\}_{b\in\B_{\bse}}$. If
$\{\psi_{\bse}\}$ and $\{\psi_{\bsv}\}$ are viewed as changing bases,
then we see that $\{\vec{w}_{\bse}\}_{{\bse}\in \hat E}$ and
$\{\vec{w}^\output_{\bse}\}_{{\bse}\in\hat E}$ determine conjugate graphs of
groups.
\end{definition}

\begin{definition}\label{d:good}
For $\bse\in\hat E$ and $\bsv=\partial_0\bse$, we say that given bases
$\B_{\bse}$ and $\B_{\bsv}$ are {\it good} if they determine
decompositions of $\group_{\bse}$ and $\group_{\bsv}$ that give a
visual simplification. Specifically, we say that $\B_{\bsv}$ and
$\B_{\bse}$ are {\it good} in any of the following four cases.
\begin{description}
\item[blowing up] There is a distinguished element
$b_{\bsv}\in\B_{\bsv}$ so that $\varphi_{\bseb}(\B_{\bseb})\subset
\langle\hat b_{\bsv}\rangle$, $\bseb\in\hat E(\bsv)$. There is no condition on
$\B_{\bse}$ in this case. (We use the notation $\hat
b_{\bsv}=\B_{\bsv}\setminus\{b_{\bsv}\}$.)
\item[unpulling] There are distinguished elements
$b_{\bsv}\in\B_{\bsv}$ and $b_{\bse}\in\B_{\bse}$ so that
$\varphi_{\bse}(\hat b_{\bse})\subset \langle \hat b_{\bsv}\rangle$,
$\varphi_{\bse}(b_{\bse})=b_{\bsv}$, and
$\varphi_{\bseb}(\group_{\bseb})\subset\langle \hat b_{\bsv}\rangle$,
$\bseb\in\hat E(\bsv)\setminus\{\bse\}$.
\item[unkilling] There is a distinguished element
$b_{\bsv}\in\B_{\bsv}$ and a partition
$\B_{\bse}=\B_{\bse}'\sqcup\B_{\bse}''$ such that
$\varphi_{\bse}(\B_{\bse}')\subset\langle\hat b_{\bsv}\rangle$,
$\varphi(\B_{\bse}'')\subset b_{\bsv}\langle\hat b_{\bsv}\rangle
b_{\bsv}^{-1}$, and $\varphi_{\bseb}(\group_{\bseb})\subset\langle\hat
b_{\bsv}\rangle$, $\bseb\in\hat E(\bsv)\setminus\{\bse\}$.
\item[\splitting] There are partitions
$\B_{\bsv}=\B_{\bsv}'\sqcup\B_{\bsv}''$ and
$\B_{\bse}=\B_{\bse}'\sqcup\B_{\bse}''$ such that
$\varphi_{\bse}(\B_{\bse}')\subset\langle\B_{\bsv}'\rangle$,
$\varphi(\B_{\bse}'')\subset\langle\B_{\bsv}''\rangle$, and for
$\bseb\in\hat E(\bsv)\setminus\{\bse\}$ either
$\varphi_{\bseb}(\group_{\bseb})\subset\langle\B_{\bsv}'\rangle$ or
$\varphi_{\bseb}(\group_{\bseb})\subset\langle\B_{\bsv}''\rangle$.
\end{description}
\end{definition}

If $\B_{\bse}$ and $\B_{\bsv}$ are good bases, then the corresponding
simplification can be performed. The bases associated to edges and
vertices of the simplified graphs are as follows. Unless explicitly
mentioned, the words $w_{{\bse},b}$ representing $\varphi_{\bse}(b)$
do not change. We use the notation as above and in
Sections~\ref{s:blowing up}--\ref{s:splitting}.
\begin{description}
\item[blowing up] After blowing up, $\B_{\bsv}=\hat b_{\bsv}$.
\item[unpulling] After unpulling, $\B_{\bse}=\hat b_{\bse}$ and
$\B_{\bsv}=\hat b_{\bsv}$.
\item[unkilling] After unkilling, $\B_{\bse'}=\B_{\bse}'$,
$\B_{\bse''}=\B_{\bse}''$, $\B_{\bsv}=\hat b_{\bsv}$, and
$w_{\bse'',b}$ represents $b_{\bsv}^{-1}\varphi_{\bse}(b) b_{\bsv}$
for $b\in\B_{\bsv}''$.
\item[\splitting] After \splitting, $\B_{\bsv'}=\B_{\bsv}'$,
$\B_{\bsv''}=\B_{\bsv}''$, $\B_{\bse'}=\B_{\bse}'$, and
$\B_{\bse''}=\B_{\bse}''$.
\end{description}

\begin{proposition}\label{p:details} 
Let $\G$ be a finite graph of finite rank free groups given as in the
beginning of this section.  Suppose that, for some $\bsv\in V$,
$\graph_{\gersten}({\conjugacyseq}(\bsv))$ can be visibly
simplified.
Then,
\begin{enumerate}
\item
 there is $\G^\output\conjugate\G$ such that $\G^\output$ can be
 simplified; and
\item
If $\vec{w}$ specifies the bonding maps of $\G$, then a conjugate
sequence $\vec{w}^\output$ may be found algorithmically so that in $\G^\output$
bases are good.
\end{enumerate}
\end{proposition}

\begin{proof}
Since (2) implies (1), it is enough to prove (2). The conjugate
sequence $\vec{w}^\output$ will be specified by supplying change of
basis automorphisms
$\{\psi_{\bseb}=\psi_{{\bseb}^{-1}}\in\Aut(\group_{\bseb})\}_{{\bseb}\in\hat
E}$ and $\{\psi_{\bsvb}\in\Aut(\group_{\bsvb})\}_{{\bsvb}\in V}$ as
well as the conjugating elements $\vec{h}=\{h_{\bseb}\in
\group_{\bseb}\}_{{\bseb}\in\hat E }$ as in
Definition~\ref{d:conjugate word sequence}. The change of basis
automorphisms can be given by specifying new bases.

Gersten's algorithm supplies $\auto\in\Aut(\group_{\bsv})$ so that
$\core(\graph_{\stallings}(\auto\groupseq(\bsv)))=
\graph_{\gersten}(\conjugacyseq(\bsv))$. By taking
$\psi_{\bsv}=\alpha$, we obtain a conjugate sequence (still denoted
$\vec{w}$) such that $\core(\graph)$ can be visibly simplified where
$\graph=\graph_{\stallings}(\groupseq(\bsv))$. We are using the
convention that all unmentioned change of basis and conjugating
automorphisms are identities.

Recall that each component of $\graph$ has a basepoint. Choose
shortest paths in each component of $\graph$ from the basepoint to the
core and take as conjugating elements the words read off along the
inverses of these paths. The resulting conjugate sequence has
basepoints in $\core(\graph)$. Hence, we may further assume that
$\graph=\core(\graph)$.

If $\graph$ can be visibly blown up, then bases are good and if
$b\in\B_{\bsv}$ is an element that does not appear as a label on
$\graph$ we take $b$ for the distinguished element $b_{\bsv}$.

In each of the remaining cases, there is a distinguished $\bse\in\hat
E(\bsv)$. Let $\graph_{\bse}$ be the component of $\graph$ indexed by
$\bse$, ie,
$\graph_{\bse}=\graph_{\stallings}(\varphi_{\bse}(\group_{\bse}))$. There
is the natural factorization
$\group_{\bse}\to\pi_1(\graph_{\bse},*)\to\group_{\bsv}$ of
$\varphi_{\bse}$ where the first map is the isomorphism coming from
Stallings algorithm (see Section~\ref{s:stallings}) and the second map
is induced by natural map of Section~\ref{s:labeled graphs}. We will
use this first map to identify $\group_{\bse}$ with
$\pi_1(\graph_{\bse},*)$. If a basis for $\pi_1(\graph_{\bse},*)$ is
given via a maximal tree $T$ for $\graph_{\bse}$, it is easy to write
this isomorphism in terms of the given bases. To make this
identification explicit, it is necessary to invert the automorphism
and this can be done algorithmically, see
Remark~\ref{r:whitehead}. Via our identification, $T$ determines a
new basis for $\group_{\bse}$ and hence a conjugate
sequence. 

If $\graph$ can be visibly unpulled, then choose an edge $\e$ of
$\graph_{\bse}$ that does not separate its component and whose
label $b\in\B_{\bsv}$ appears exactly once in $\graph$. 
Choose a maximal tree $T$ for $\graph_{\bse}$ so that $\e\not\subset
T$ and let $\B_{\bse}^\output$ be the basis determined by
$T$. There is an element $b_{\bse}\in\B_{\bse}^\output$
corresponding to $\e$. Let $b_{\bsv}$ be $\varphi_{\bse}$--image in
$\group_{\bsv}$ of $b_{\bse}$. Set
$\B_{\bsv}^\output=\{b_{\bsv}\}\sqcup\hat b$ which is a basis for
$\group_{\bsv}$ since $b$ appears exactly once in $b_{\bsv}$ when
expressed as a $\B_{\bsv}$--word. The new
bases are good and determine change of basis automorphisms giving rise
to the desired $\vec{w}^\output$.

If $\graph$ can be visibly unkilled, then choose an edge $\e$ of
$\graph_{\bse}$ that separates its component and whose label
$b\in\B_{\bsv}$ appears exactly once in $\graph$. By changing the
orientation of $\e$ and inverting $b$ if necessary, we may assume that
$\partial_0\e$ and the basepoint are in the same component of the
graph obtained by removing $\e$ from $\graph_{\bse}$. Let $T$ be a
maximal tree for $\graph_{\bse}$. This gives the desired new basis
$\B^\output_{\bse}$. Choosing a shortest path in $T$ from the
$\partial_0\e$ to the basepoint gives rise to a conjugating element
$h_{\bse}$ which has the effect of changing the basepoint of
$\graph_{\bse}$ to $\partial_0\e$. The new bases are good. Indeed, the
partition of $\B^\output_{\bse}$ is induced by the separating
edge. More precisely, let $(\B_{\bse}^\output)''$ be the
elements of $\B^\output_{\bse}$ that contain the letter $b$ and let
$(\B_{\bse}^{\output})'$ be the complement. Set $b_{\bsv}=b$.

If $\graph$ can be visibly \splited, then write
$\graph=\graph'\vee\graph''$ respecting the non-trivial partition
$\B_{\bsv}'\sqcup\B_{\bsv}''$ of $\B_{\bsv}$.  A choice of maximal
tree $T$ for $\graph_{\bse}$ gives the desired new basis
$\B_{\bse}^\output$ and a shortest path from the wedge point to
the basepoint determines a conjugating element $h_{\bse}$. With these
choices, the new bases are good. Indeed, the partition of
$\B^\output_{\bse}$ is induced by the wedge. More precisely,
$(\B_{\bse}^{\output})'$ corresponds to the set of edges of
$\graph'\cap\graph_{\bse}$ not in $T$ and
$(\B_{\bse}^{\output})''$ corresponds to the remaining
edges not in $T$.
\end{proof}

\begin{example}
Let $\G$ have underlying graph $\bsgraph$ as in
\figref{f:algorithm.example}. Suppose that
$\B_{\bse}=\{a_1,a_2 \}$, $\B_{\bsv}=\{b_1,b_2 \}$, and
$\B_f=\{ z\}$. We will not specify $G_\bsvb$ since it will not
change. Suppose that $\varphi_{\bse}(a_1)=b_1^2b_2^2$,
$\varphi_{\bse}(a_2)=b_1^2b_2^2b_1^2$, $\varphi_{\bseb}(z)=b_1$,
$\varphi_{\bseb^{-1}}(z)=b_2$. Then,
$\graph(\bsv)=\graph_{\stallings}(\groupseq(\bsv))$ is a Gersten
representative, ie, $\graph_{\stallings}(\groupseq(\bsv))=
\graph_{\gersten}(\conjugacyseq(\bsv))$.  $\graph(\bsv)$ is
displayed in \figref{f:algorithm.example}, and it can be visibly
\splited. The given bases are not good. A good basis for $\group_{\bse}$
corresponds to the one determined by the wedge point in
$\graph(\bsv)$. The change of basis automorphism
$\psi_{\bse}\in\Aut(\group_{\bse})$ is given by
$\psi_{\bse}(a_1)=a_1^{-1}a_2$ and $\psi_{\bse}(a_2)=a_2^{-1}a_1a_1$, ie,
$\varphi_{\bse}\circ\psi_{\bse}(a_1)=b_1^2$ and
$\varphi_{\bse}\circ\psi_{\bse}(a_2)=b_2^2$.

If $\G'$ is the result of \splitting\ $\G$, then $\B_{\bse'}=\{a_1\}$,
$\B_{\bse''}=\{a_2\}$, $\B_{\bsv'}=\{b_1\}$, $\B_{\bsv''}=\{b_2\}$,
$\varphi_{\bse'}(a_1)=b_1^2$, $\varphi_{\bse''}(a_2)=b_2^2$,
$\varphi_{(\bse')^{-1}}(a_1)=\varphi_{\bse^{-1}}\circ\psi_{\bse}(a_1)$, and
$\varphi_{(\bse'')^{-1}}(a_2)=\varphi_{\bse^{-1}}\circ\psi_{\bse}(a_2)$.

The next step in the algorithm would be to reduce $\G'$. The example
could have been complicated by post-composing $\varphi_{\bse}$,
$\varphi_{\bseb}$, and $\varphi_{\bseb^{-1}}$ by some
$\psi_{\bsv}\in\Aut(\group_{\bsv})$. In that case, we would use Gersten's
algorithm first (and discover $\psi_{\bsv}$).
\end{example}

\begin{figure}[ht!]\anchor{f:algorithm.example}
\cl{\scalebox{0.65}{\input{e1.pstex_t}}}
\cl{\scalebox{0.65}{\input{e1.1.pstex_t}}}
\nocolon\caption{}\label{f:algorithm.example}
\end{figure}

\end{document}

%% file: gtoutput.tex
\def\ifplaintex{\expandafter\ifx\csname documentclass\endcsname\relax}

\ifplaintex 
\else
\fi

\expandafter\ifx\csname beginpicture\endcsname\relax
\expandafter\ifx\csname documentclass\endcsname\relax
\input pictex \else
\input prepictex \input pictex \input postpictex \fi\fi

\def\gt{{\mathsurround=0pt\it $\cal G\mskip-2mu$eometry \&\ 
$\cal T\!\!$opology}}        

\def\gtp{{\mathsurround=0pt\it $\cal G\mskip-2mu$eometry \&\ 
$\cal T\!\!$opology $\cal P\!$ublications}}  

\def\lognumber#1{\def\thelognumber{#1}}
\def\volumenumber#1{\def\thevolumenumber{#1}}
\def\papernumber#1{\def\thepapernumber{#1}}
\def\volumeyear#1{\def\thevolumeyear{#1}}

\def\pagenumbers#1#2{\def\startpage{#1}\def\finishpage{#2}}
\def\published#1{\def\publishdate{#1}}
\def\proposed#1{\def\theproposer{#1}}
\def\seconded#1{\def\theseconders{#1}}
\def\received#1{\def\receiveddate{#1}}
\def\revised#1{\def\reviseddate{#1}}
\def\accepted#1{\def\accepteddate{#1}}

\def\asciiaddress#1{\def\theasciiaddress{#1}}
\def\asciiemail#1{\def\theasciiemail{#1}}

\long\def\asciiabstract#1{\long\def\theasciiabstract{#1}}

\let\\\par\let\thelognumber\relax
\let\thevolumenumber\relax\let\thepapernumber\relax
\let\thevolumeyear\relax\let\thesamplenumber\relax\let\startpage\relax
\let\finishpage\relax\let\publishdate\relax\let\receiveddate\relax
\let\reviseddate\relax\let\accepteddate\relax\let\theasciititle\relax
\let\theasciiauthors\relax\let\theasciiaddress\relax
\let\theasciiabstract\relax
\let\theasciiemail\relax\let\theshortauthors\relax\let\theshorttitle\relax

\long\def\maketitlep{   

\count0=\startpage

\gt\hfill      
\beginpicture
\setcoordinatesystem units <0.33truein, 0.33truein> point at 2.2 0.9
\setplotsymbol ({$\cal G$})
\plotsymbolspacing=9truept
\circulararc 315 degrees from 0 1 center at 0 0
\setplotsymbol ({$\cal T$})
\circulararc 315 degrees from 1 -1 center at 1 0
\endpicture
\break
{\small\ifx\thesamplenumber\relax 
Volume \else Sample
\fi\thevolumenumber\ (\thevolumeyear)
\startpage--\finishpage\nl
Published: \publishdate}
\vglue 0.5truein plus 0.4fil minus 0.1truein

{\parskip=0pt\leftskip 0pt plus 1fil\def\\{\par\smallskip}{\ifplaintex\large
\else\Large\fi\bf\thetitle}\par\medskip}   

\vglue 0pt plus 0.1fil 

{\parskip=0pt\leftskip 0pt plus 1fil\def\\{\par}{\sc\theauthors}
\par\medskip}

\vglue 0pt plus 0.1fil 

{\small\parskip=0pt\let\newline\\
{\leftskip 0pt plus 1fil\def\\{\par}{\sl\theaddress}\par}
\expandafter\ifx\theemail\relax    
\relax\else\vglue 5pt plus 0.02fil minus 2pt\def\\{\stdspace{\rm 
and}\stdspace} 
\cl{Email:\stdspace\tt\theemail}\fi
\ifx\theurl\relax                  
\relax\else\vglue 5pt plus 0.02fil minus 2pt\def\\{\stdspace{\rm 
and}\stdspace}
\cl{URL:\stdspace\tt\theurl}\fi\par}

\vglue 7pt plus 0.3fil minus 3pt

{\bf Abstract}
\vglue 5pt plus 0.1fil minus 2pt

\theabstract

\vglue 7pt plus 0.3fil minus 3pt

{\bf AMS Classification numbers}\quad Primary:\quad \theprimaryclass

Secondary:\quad \thesecondaryclass

\vglue 5pt plus 0.3fil minus 2pt

{\bf Keywords:}\quad \thekeywords

\vglue 10pt plus 0.5fil minus 5pt

{\small  Proposed: \theproposer\hfill Received: \receiveddate\nl
Seconded: \theseconders\hfill 
\ifx\reviseddate\relax                         
Accepted: \accepteddate                        
\else
Revised: \reviseddate                          
\fi}
\eject
}       

\let\maketitlepage\maketitlep
\let\maketitle\maketitlepage

\font\phead=cmsl9 scaled 950
\font=cmsl9 scaled 1050
\font\pnum=cmbx10 scaled 913
\font\lnum=cmbx10 
\font\pfoot=cmsl9 scaled 950
\font=cmsl9 scaled 1050
\ifplaintex
\headline{\vbox to 0pt{\vskip -4.5mm\line{\small\phead\ifnum
\count0=\startpage ISSN 1364-0380 (on line)
1465-3060 (printed) \hfill {\pnum\folio}\else\ifodd\count0\def\\{ }%
\ifx\theshorttitle\relax\thetitle\else\theshorttitle\fi\hfill{\pnum\folio}
\else\def\\{ and }{\pnum\folio}\hfill\ifx\theshortauthors\relax\theauthors
\else\theshortauthors\fi\fi\fi}\vss}}
\footline{\vbox to 0pt{\vglue 0mm\line{\small\pfoot\ifnum\count0=\startpage
\copyright\ \gtp\hfill\else
\gt, Volume \thevolumenumber\ (\thevolumeyear)\hfill\fi}\vss
}}
\else
\makeatletter
\def\@oddhead{{\small\lhead\ifnum\count0=\startpage ISSN 1364-0380 (on line)
1465-3060 (printed) \hfill {\lnum\number\count0}\else\ifodd\count0
\def\\{ }\ifx\theshorttitle\relax \thetitle \else\theshorttitle\fi\hfill
{\lnum\number\count0}\else\def\\{ and }{\lnum\number\count0}
\hfill\ifx\theshortauthors\relax 
\theauthors\else\theshortauthors\fi\fi\fi}}\def\@evenhead{\@oddhead}
\def\@oddfoot{\small\lfoot\ifnum\count0=\startpage\copyright\ \gtp\hfill\else
\gt, Volume \thevolumenumber\ (\thevolumeyear)\hfill\fi}
\def\@evenfoot{\@oddfoot}
\makeatother
\fi

\newwrite\gtoutfile
\long\gdef\makeheadfile{  
{\def\\{, }\def\s{ }
\immediate\openout\gtoutfile head.xxx
\immediate\write\gtoutfile{Proxy-for: \ifx\theasciiauthors\relax
\theauthors\else\theasciiauthors\fi\s<\ifx\theasciiemail\relax\theemail\else\theasciiemail\fi>}
\immediate\write\gtoutfile{\noexpand\\}
\immediate\write\gtoutfile{Authors: \ifx\theasciiauthors\relax
\theauthors\else\theasciiauthors\fi}
{\def\\{ }\immediate\write\gtoutfile{Title: \ifx\theasciititle\relax
\thetitle\else\theasciititle\fi}}
\immediate\write\gtoutfile{Subj-class: GT or SG or MG etc}
\immediate\write\gtoutfile{MSC-class: \theprimaryclass\ifx\thesecondaryclass\relax\else, \thesecondaryclass\fi}
\immediate\write\gtoutfile{Journal-ref: Geom. Topol. \thevolumenumber
(\thevolumeyear) \startpage-\finishpage}
\immediate\write\gtoutfile{Comments: Published by Geometry and Topology at}
\immediate\write\gtoutfile{\s\s http://www.maths.warwick.ac.uk/gt/GTVol\thevolumenumber/paper\thepapernumber.abs.html}
\immediate\write\gtoutfile{\noexpand\\}
\immediate\write\gtoutfile{}
\ifx\theasciiabstract\relax
\immediate\write\gtoutfile{\theabstract}\else
\immediate\write\gtoutfile{\theasciiabstract}\fi
\immediate\write\gtoutfile{}
\immediate\write\gtoutfile{\noexpand\\}
\immediate\write\gtoutfile{}
\immediate\closeout\gtoutfile}}  

\def\maketitlepage{\maketitlep\makeheadfile}
\let\maketitle\maketitlepage

%% file: reduce.pstex_t
\begin{picture}(0,0)%
\includegraphics{reduce.pstex}%
\end{picture}%
\setlength{\unitlength}{3947sp}%
\begingroup\makeatletter\ifx\SetFigFont\undefined%
\gdef\SetFigFont#1#2#3#4#5{%
  \reset@font\fontsize{#1}{#2pt}%
  \fontfamily{#3}\fontseries{#4}\fontshape{#5}%
  \selectfont}%
\fi\endgroup%
\begin{picture}(4899,2325)(64,-2374)
\put(2026,-2011){\makebox(0,0)[b]{\smash{{\SetFigFont{14}{16.8}{\rmdefault}{\mddefault}{\updefault}{\color[rgb]{1,0,0}$\group_{\bse}$}%
}}}}
\put(1351,-2311){\makebox(0,0)[b]{\smash{{\SetFigFont{14}{16.8}{\rmdefault}{\mddefault}{\updefault}{\color[rgb]{1,0,0}$\group_{\bse}$}%
}}}}
\put(2026,-1186){\makebox(0,0)[b]{\smash{{\SetFigFont{14}{16.8}{\rmdefault}{\mddefault}{\updefault}{\color[rgb]{1,0,0}$\group_{\bse}$}%
}}}}
\put(751,-2011){\makebox(0,0)[b]{\smash{{\SetFigFont{14}{16.8}{\rmdefault}{\mddefault}{\updefault}$\group_{\bseb}$}}}}
\put(4351,-2011){\makebox(0,0)[b]{\smash{{\SetFigFont{14}{16.8}{\rmdefault}{\mddefault}{\updefault}$\group_{\bseb}$}}}}
\put(2701,-811){\makebox(0,0)[b]{\smash{{\SetFigFont{14}{16.8}{\rmdefault}{\mddefault}{\updefault}{\color[rgb]{1,0,0}$\group_{\bse}$}%
}}}}
\end{picture}%

%% file: blowup2.pstex_t
\begin{picture}(0,0)%
\includegraphics{blowup2.pstex}%
\end{picture}%
\setlength{\unitlength}{3947sp}%
\begingroup\makeatletter\ifx\SetFigFont\undefined%
\gdef\SetFigFont#1#2#3#4#5{%
  \reset@font\fontsize{#1}{#2pt}%
  \fontfamily{#3}\fontseries{#4}\fontshape{#5}%
  \selectfont}%
\fi\endgroup%
\begin{picture}(5425,1226)(589,-1574)
\put(5776,-1036){\makebox(0,0)[b]{\smash{{\SetFigFont{14}{16.8}{\rmdefault}{\mddefault}{\updefault}{\color[rgb]{1,0,0}$\trivialgroup$}%
}}}}
\put(5101,-886){\makebox(0,0)[b]{\smash{{\SetFigFont{14}{16.8}{\rmdefault}{\mddefault}{\updefault}{\color[rgb]{1,0,0}$\group_{\bsv}'$}%
}}}}
\put(1876,-736){\makebox(0,0)[b]{\smash{{\SetFigFont{14}{16.8}{\rmdefault}{\mddefault}{\updefault}{\color[rgb]{1,0,0}$\group_{\bsv}'*\langle t\rangle$}%
}}}}
\end{picture}%

%% file: blowup3.pstex_t
\begin{picture}(0,0)%
\includegraphics{blowup3.pstex}%
\end{picture}%
\setlength{\unitlength}{3947sp}%
\begingroup\makeatletter\ifx\SetFigFont\undefined%
\gdef\SetFigFont#1#2#3#4#5{%
  \reset@font\fontsize{#1}{#2pt}%
  \fontfamily{#3}\fontseries{#4}\fontshape{#5}%
  \selectfont}%
\fi\endgroup%
\begin{picture}(4062,1974)(589,-2023)
\put(4576,-586){\makebox(0,0)[b]{\smash{{\SetFigFont{14}{16.8}{\rmdefault}{\mddefault}{\updefault}{\color[rgb]{1,0,0}$\group_{\bsv}''$}%
}}}}
\put(4576,-1636){\makebox(0,0)[b]{\smash{{\SetFigFont{14}{16.8}{\rmdefault}{\mddefault}{\updefault}{\color[rgb]{1,0,0}$\group_{\bsv}'$}%
}}}}
\put(4651,-1036){\makebox(0,0)[b]{\smash{{\SetFigFont{14}{16.8}{\rmdefault}{\mddefault}{\updefault}{\color[rgb]{1,0,0}$\trivialgroup$}%
}}}}
\put(1876,-736){\makebox(0,0)[b]{\smash{{\SetFigFont{14}{16.8}{\rmdefault}{\mddefault}{\updefault}{\color[rgb]{1,0,0}$\group_{\bsv}'*\group_\bsv''$}%
}}}}
\end{picture}%

%% file: unpull.pstex_t
\begin{picture}(0,0)%
\includegraphics{unpull.pstex}%
\end{picture}%
\setlength{\unitlength}{3947sp}%
\begingroup\makeatletter\ifx\SetFigFont\undefined%
\gdef\SetFigFont#1#2#3#4#5{%
  \reset@font\fontsize{#1}{#2pt}%
  \fontfamily{#3}\fontseries{#4}\fontshape{#5}%
  \selectfont}%
\fi\endgroup%
\begin{picture}(5874,1449)(439,-1723)
\put(5701,-1186){\makebox(0,0)[b]{\smash{{\SetFigFont{14}{16.8}{\rmdefault}{\mddefault}{\updefault}{\color[rgb]{1,0,0}$\group_{\bse}'$}%
}}}}
\put(5176,-736){\makebox(0,0)[b]{\smash{{\SetFigFont{14}{16.8}{\rmdefault}{\mddefault}{\updefault}{\color[rgb]{1,0,0}$\group_{\bsv}'$}%
}}}}
\put(1876,-736){\makebox(0,0)[b]{\smash{{\SetFigFont{14}{16.8}{\rmdefault}{\mddefault}{\updefault}{\color[rgb]{1,0,0}$\group_{\bsv}'*\mathbb Z$}%
}}}}
\put(2401,-1186){\makebox(0,0)[b]{\smash{{\SetFigFont{14}{16.8}{\rmdefault}{\mddefault}{\updefault}{\color[rgb]{1,0,0}$\group_{\bse}'*\mathbb Z$}%
}}}}
\end{picture}%

%% file: unkill.pstex_t
\begin{picture}(0,0)%
\includegraphics{unkill.pstex}%
\end{picture}%
\setlength{\unitlength}{3947sp}%
\begingroup\makeatletter\ifx\SetFigFont\undefined%
\gdef\SetFigFont#1#2#3#4#5{%
  \reset@font\fontsize{#1}{#2pt}%
  \fontfamily{#3}\fontseries{#4}\fontshape{#5}%
  \selectfont}%
\fi\endgroup%
\begin{picture}(7449,1524)(64,-1723)
\put(5551,-511){\makebox(0,0)[b]{\smash{{\SetFigFont{14}{16.8}{\rmdefault}{\mddefault}{\updefault}{\color[rgb]{1,0,0}$\group_{\bse}'$}%
}}}}
\put(2026,-1186){\makebox(0,0)[b]{\smash{{\SetFigFont{14}{16.8}{\rmdefault}{\mddefault}{\updefault}{\color[rgb]{1,0,0}$\group_{\bse}'*\group_{\bse}''$}%
}}}}
\put(1651,-736){\makebox(0,0)[b]{\smash{{\SetFigFont{14}{16.8}{\rmdefault}{\mddefault}{\updefault}{\color[rgb]{1,0,0}$\group_{\bsv}'*\langle t \rangle$}%
}}}}
\put(4951,-736){\makebox(0,0)[b]{\smash{{\SetFigFont{14}{16.8}{\rmdefault}{\mddefault}{\updefault}{\color[rgb]{1,0,0}$\group_{\bsv}'$}%
}}}}
\put(5551,-1486){\makebox(0,0)[b]{\smash{{\SetFigFont{14}{16.8}{\rmdefault}{\mddefault}{\updefault}{\color[rgb]{1,0,0}$\group_{\bse}''$}%
}}}}
\end{picture}%

%% file: split.pstex_t
\begin{picture}(0,0)%
\includegraphics{split.pstex}%
\end{picture}%
\setlength{\unitlength}{3947sp}%
\begingroup\makeatletter\ifx\SetFigFont\undefined%
\gdef\SetFigFont#1#2#3#4#5{%
  \reset@font\fontsize{#1}{#2pt}%
  \fontfamily{#3}\fontseries{#4}\fontshape{#5}%
  \selectfont}%
\fi\endgroup%
\begin{picture}(6988,2049)(64,-2098)
\put(4951,-586){\makebox(0,0)[b]{\smash{{\SetFigFont{14}{16.8}{\rmdefault}{\mddefault}{\updefault}{\color[rgb]{1,0,0}$\group_{\bsv}'$}%
}}}}
\put(5551,-811){\makebox(0,0)[b]{\smash{{\SetFigFont{14}{16.8}{\rmdefault}{\mddefault}{\updefault}{\color[rgb]{1,0,0}$\group_{\bse}'$}%
}}}}
\put(2026,-1186){\makebox(0,0)[b]{\smash{{\SetFigFont{14}{16.8}{\rmdefault}{\mddefault}{\updefault}{\color[rgb]{1,0,0}$\group_{\bse}'*\group_{\bse}''$}%
}}}}
\put(4951,-1636){\makebox(0,0)[b]{\smash{{\SetFigFont{14}{16.8}{\rmdefault}{\mddefault}{\updefault}{\color[rgb]{1,0,0}$\group_{\bsv}''$}%
}}}}
\put(5551,-1411){\makebox(0,0)[b]{\smash{{\SetFigFont{14}{16.8}{\rmdefault}{\mddefault}{\updefault}{\color[rgb]{1,0,0}$\group_{\bse}''$}%
}}}}
\put(1651,-736){\makebox(0,0)[b]{\smash{{\SetFigFont{14}{16.8}{\rmdefault}{\mddefault}{\updefault}{\color[rgb]{1,0,0}$\group_{\bsv}'*\group_{\bsv}''$}%
}}}}
\end{picture}%

%% file: gersten.pstex_t
\begin{picture}(0,0)%
\includegraphics{gersten.pstex}%
\end{picture}%
\setlength{\unitlength}{3947sp}%
\begingroup\makeatletter\ifx\SetFigFont\undefined%
\gdef\SetFigFont#1#2#3#4#5{%
  \reset@font\fontsize{#1}{#2pt}%
  \fontfamily{#3}\fontseries{#4}\fontshape{#5}%
  \selectfont}%
\fi\endgroup%
\begin{picture}(3624,4805)(499,-4269)
\put(3931,-2506){\makebox(0,0)[b]{\smash{{\SetFigFont{12}{14.4}{\rmdefault}{\mddefault}{\updefault}{\color[rgb]{0,0,0}$(\graph_{\stallings},*)$}%
}}}}
\put(4096,-4141){\makebox(0,0)[b]{\smash{{\SetFigFont{12}{14.4}{\rmdefault}{\mddefault}{\updefault}{\color[rgb]{0,0,0}$\graph_{\gersten}$}%
}}}}
\put(3916,-1276){\makebox(0,0)[b]{\smash{{\SetFigFont{12}{14.4}{\rmdefault}{\mddefault}{\updefault}{\color[rgb]{0,0,0}$(\graph,*)$}%
}}}}
\end{picture}%

%% file: flow2.pstex_t
\begin{picture}(0,0)%
\includegraphics{flow2.pstex}%
\end{picture}%
\setlength{\unitlength}{3947sp}%
\begingroup\makeatletter\ifx\SetFigFont\undefined%
\gdef\SetFigFont#1#2#3#4#5{%
  \reset@font\fontsize{#1}{#2pt}%
  \fontfamily{#3}\fontseries{#4}\fontshape{#5}%
  \selectfont}%
\fi\endgroup%
\begin{picture}(5045,4070)(1789,-3748)
\put(3466,-1546){\makebox(0,0)[b]{\smash{{\SetFigFont{12}{14.4}{\rmdefault}{\mddefault}{\updefault}{\color[rgb]{0,0,0}$\exists$? $\bsv$}%
}}}}
\put(3451,-136){\makebox(0,0)[b]{\smash{{\SetFigFont{12}{14.4}{\rmdefault}{\mddefault}{\updefault}{\color[rgb]{0,0,0}Reduce $\G$}%
}}}}
\put(5776,-361){\makebox(0,0)[b]{\smash{{\SetFigFont{12}{14.4}{\rmdefault}{\mddefault}{\updefault}{\color[rgb]{0,0,0}Input $\G$}%
}}}}
\put(3481,-3391){\makebox(0,0)[b]{\smash{{\SetFigFont{12}{14.4}{\rmdefault}{\mddefault}{\updefault}{\color[rgb]{0,0,0}Simplify}%
}}}}
\put(3421,-2026){\makebox(0,0)[b]{\smash{{\SetFigFont{12}{14.4}{\rmdefault}{\mddefault}{\updefault}{\color[rgb]{0,0,0}simplifies}%
}}}}
\put(3451,-1786){\makebox(0,0)[b]{\smash{{\SetFigFont{12}{14.4}{\rmdefault}{\mddefault}{\updefault}{\color[rgb]{0,0,0}s.t. $\graph_{\gersten}(\conjugacyseq(\bsv))$ visibly}%
}}}}
\put(5926,-1786){\makebox(0,0)[b]{\smash{{\SetFigFont{12}{14.4}{\rmdefault}{\mddefault}{\updefault}{\color[rgb]{0,0,0}DONE}%
}}}}
\end{picture}%

%% file: model.pstex_t
\begin{picture}(0,0)%
\includegraphics{model.pstex}%
\end{picture}%
\setlength{\unitlength}{3947sp}%
\begingroup\makeatletter\ifx\SetFigFont\undefined%
\gdef\SetFigFont#1#2#3#4#5{%
  \reset@font\fontsize{#1}{#2pt}%
  \fontfamily{#3}\fontseries{#4}\fontshape{#5}%
  \selectfont}%
\fi\endgroup%
\begin{picture}(4800,2414)(901,-2768)
\put(3826,-1636){\makebox(0,0)[b]{\smash{\SetFigFont{17}{20.4}{\rmdefault}{\mddefault}{\updefault}$\longrightarrow$}}}
\put(901,-1636){\makebox(0,0)[b]{\smash{\SetFigFont{17}{20.4}{\rmdefault}{\mddefault}{\updefault}$\bar X$}}}
\put(5701,-1636){\makebox(0,0)[lb]{\smash{\SetFigFont{17}{20.4}{\rmdefault}{\mddefault}{\updefault}$\bar S$}}}
\end{picture}

%% file: 0collapse.pstex_t
\begin{picture}(0,0)%
\includegraphics{0collapse.pstex}%
\end{picture}%
\setlength{\unitlength}{3947sp}%
\begingroup\makeatletter\ifx\SetFigFont\undefined%
\gdef\SetFigFont#1#2#3#4#5{%
  \reset@font\fontsize{#1}{#2pt}%
  \fontfamily{#3}\fontseries{#4}\fontshape{#5}%
  \selectfont}%
\fi\endgroup%
\begin{picture}(4512,1224)(901,-973)
\put(1201,-211){\makebox(0,0)[b]{\smash{\SetFigFont{17}{20.4}{\rmdefault}{\mddefault}{\updefault}{\color[rgb]{0,0,0}$C_0$}%
}}}
\put(901,-886){\makebox(0,0)[b]{\smash{\SetFigFont{17}{20.4}{\rmdefault}{\mddefault}{\updefault}{\color[rgb]{0,0,0}$X$}%
}}}
\put(4351,-811){\makebox(0,0)[b]{\smash{\SetFigFont{17}{20.4}{\rmdefault}{\mddefault}{\updefault}{\color[rgb]{0,0,0}$\result(X)$}%
}}}
\put(3751,-61){\makebox(0,0)[b]{\smash{\SetFigFont{12}{14.4}{\rmdefault}{\mddefault}{\updefault}{\color[rgb]{0,0,0}0-Simplify}%
}}}
\put(1801,-661){\makebox(0,0)[b]{\smash{\SetFigFont{17}{20.4}{\rmdefault}{\mddefault}{\updefault}{\color[rgb]{0,0,0}$C$}%
}}}
\end{picture}

%% file: 1collapse.pstex_t
\begin{picture}(0,0)%
\includegraphics{1collapse.pstex}%
\end{picture}%
\setlength{\unitlength}{3947sp}%
\begingroup\makeatletter\ifx\SetFigFont\undefined%
\gdef\SetFigFont#1#2#3#4#5{%
  \reset@font\fontsize{#1}{#2pt}%
  \fontfamily{#3}\fontseries{#4}\fontshape{#5}%
  \selectfont}%
\fi\endgroup%
\begin{picture}(4212,2499)(976,-2248)
\put(976,-1036){\makebox(0,0)[b]{\smash{\SetFigFont{17}{20.4}{\rmdefault}{\mddefault}{\updefault}{\color[rgb]{0,0,0}$C_0$}%
}}}
\put(3826,-811){\makebox(0,0)[b]{\smash{\SetFigFont{12}{14.4}{\rmdefault}{\mddefault}{\updefault}{\color[rgb]{0,0,0}I-Simplify}%
}}}
\put(976,-286){\makebox(0,0)[b]{\smash{\SetFigFont{17}{20.4}{\rmdefault}{\mddefault}{\updefault}{\color[rgb]{0,0,0}$X$}%
}}}
\put(4426,-286){\makebox(0,0)[b]{\smash{\SetFigFont{17}{20.4}{\rmdefault}{\mddefault}{\updefault}{\color[rgb]{0,0,0}$\result(X)$}%
}}}
\put(1801,-1036){\makebox(0,0)[b]{\smash{\SetFigFont{17}{20.4}{\rmdefault}{\mddefault}{\updefault}{\color[rgb]{0,0,0}$C$}%
}}}
\end{picture}

%% file: 2collapse.pstex_t
\begin{picture}(0,0)%
\includegraphics{2collapse.pstex}%
\end{picture}%
\setlength{\unitlength}{3947sp}%
\begingroup\makeatletter\ifx\SetFigFont\undefined%
\gdef\SetFigFont#1#2#3#4#5{%
  \reset@font\fontsize{#1}{#2pt}%
  \fontfamily{#3}\fontseries{#4}\fontshape{#5}%
  \selectfont}%
\fi\endgroup%
\begin{picture}(6024,2424)(589,-2173)
\put(3901,-811){\makebox(0,0)[b]{\smash{\SetFigFont{12}{14.4}{\rmdefault}{\mddefault}{\updefault}{\color[rgb]{0,0,0}II-Simplify}%
}}}
\put(976,-1036){\makebox(0,0)[b]{\smash{\SetFigFont{17}{20.4}{\rmdefault}{\mddefault}{\updefault}{\color[rgb]{0,0,0}$C_0$}%
}}}
\put(601,-361){\makebox(0,0)[b]{\smash{\SetFigFont{17}{20.4}{\rmdefault}{\mddefault}{\updefault}{\color[rgb]{0,0,0}$X$}%
}}}
\put(4351,-361){\makebox(0,0)[b]{\smash{\SetFigFont{17}{20.4}{\rmdefault}{\mddefault}{\updefault}{\color[rgb]{0,0,0}$\result(X)$}%
}}}
\put(1801,-1036){\makebox(0,0)[b]{\smash{\SetFigFont{17}{20.4}{\rmdefault}{\mddefault}{\updefault}{\color[rgb]{0,0,0}$C$}%
}}}
\end{picture}

%% file: 3collapse.pstex_t
\begin{picture}(0,0)%
\includegraphics{3collapse.pstex}%
\end{picture}%
\setlength{\unitlength}{3947sp}%
\begingroup\makeatletter\ifx\SetFigFont\undefined%
\gdef\SetFigFont#1#2#3#4#5{%
  \reset@font\fontsize{#1}{#2pt}%
  \fontfamily{#3}\fontseries{#4}\fontshape{#5}%
  \selectfont}%
\fi\endgroup%
\begin{picture}(6024,2424)(589,-2173)
\put(3901,-886){\makebox(0,0)[b]{\smash{\SetFigFont{12}{14.4}{\rmdefault}{\mddefault}{\updefault}{\color[rgb]{0,0,0}III-Simplify}%
}}}
\put(4276,-361){\makebox(0,0)[b]{\smash{\SetFigFont{17}{20.4}{\rmdefault}{\mddefault}{\updefault}{\color[rgb]{0,0,0}$\result(X)$}%
}}}
\put(976,-1036){\makebox(0,0)[b]{\smash{\SetFigFont{17}{20.4}{\rmdefault}{\mddefault}{\updefault}{\color[rgb]{0,0,0}$C_0$}%
}}}
\put(601,-361){\makebox(0,0)[b]{\smash{\SetFigFont{17}{20.4}{\rmdefault}{\mddefault}{\updefault}{\color[rgb]{0,0,0}$X$}%
}}}
\put(1801,-1036){\makebox(0,0)[b]{\smash{\SetFigFont{17}{20.4}{\rmdefault}{\mddefault}{\updefault}{\color[rgb]{0,0,0}$C$}%
}}}
\end{picture}

%% file: blowup.pstex_t
\begin{picture}(0,0)%
\includegraphics{blowup.pstex}%
\end{picture}%
\setlength{\unitlength}{3947sp}%
\begingroup\makeatletter\ifx\SetFigFont\undefined%
\gdef\SetFigFont#1#2#3#4#5{%
  \reset@font\fontsize{#1}{#2pt}%
  \fontfamily{#3}\fontseries{#4}\fontshape{#5}%
  \selectfont}%
\fi\endgroup%
\begin{picture}(5862,4224)(751,-3973)
\put(3301,-2086){\makebox(0,0)[b]{\smash{\SetFigFont{12}{14.4}{\rmdefault}{\mddefault}{\updefault}{\color[rgb]{0,0,0}Blow up}%
}}}
\put(1576,-1036){\makebox(0,0)[b]{\smash{\SetFigFont{17}{20.4}{\rmdefault}{\mddefault}{\updefault}{\color[rgb]{0,0,0}$C$}%
}}}
\put(5401,-886){\makebox(0,0)[b]{\smash{\SetFigFont{17}{20.4}{\rmdefault}{\mddefault}{\updefault}{\color[rgb]{0,0,0}$C$}%
}}}
\put(826,-1036){\makebox(0,0)[b]{\smash{\SetFigFont{17}{20.4}{\rmdefault}{\mddefault}{\updefault}{\color[rgb]{0,0,0}$X$}%
}}}
\put(751,-3211){\makebox(0,0)[b]{\smash{\SetFigFont{17}{20.4}{\rmdefault}{\mddefault}{\updefault}{\color[rgb]{0,0,0}$S$}%
}}}
\put(3751,-1036){\makebox(0,0)[b]{\smash{\SetFigFont{17}{20.4}{\rmdefault}{\mddefault}{\updefault}{\color[rgb]{0,0,0}$\result(X)$}%
}}}
\put(3751,-3211){\makebox(0,0)[b]{\smash{\SetFigFont{17}{20.4}{\rmdefault}{\mddefault}{\updefault}{\color[rgb]{0,0,0}$\result(S)$}%
}}}
\end{picture}

%% file: e1.pstex_t
\begin{picture}(0,0)%
\includegraphics{e1.pstex}%
\end{picture}%
\setlength{\unitlength}{3947sp}%
\begingroup\makeatletter\ifx\SetFigFont\undefined%
\gdef\SetFigFont#1#2#3#4#5{%
  \reset@font\fontsize{#1}{#2pt}%
  \fontfamily{#3}\fontseries{#4}\fontshape{#5}%
  \selectfont}%
\fi\endgroup%
\begin{picture}(5821,1971)(555,-1990)
\put(3151,-886){\makebox(0,0)[b]{\smash{{\SetFigFont{14}{16.8}{\rmdefault}{\mddefault}{\updefault}{\color[rgb]{0,0,0}$\bseb$}%
}}}}
\put(826,-211){\makebox(0,0)[b]{\smash{{\SetFigFont{17}{20.4}{\rmdefault}{\mddefault}{\updefault}{\color[rgb]{0,0,0}$\G$}%
}}}}
\put(4276,-211){\makebox(0,0)[b]{\smash{{\SetFigFont{17}{20.4}{\rmdefault}{\mddefault}{\updefault}{\color[rgb]{0,0,0}$\G'$}%
}}}}
\put(4426,-886){\makebox(0,0)[b]{\smash{{\SetFigFont{14}{16.8}{\rmdefault}{\mddefault}{\updefault}{\color[rgb]{0,0,0}$\bsvb$}%
}}}}
\put(4951,-1636){\makebox(0,0)[b]{\smash{{\SetFigFont{14}{16.8}{\rmdefault}{\mddefault}{\updefault}{\color[rgb]{0,0,0}$\bse''$}%
}}}}
\put(4951,-586){\makebox(0,0)[b]{\smash{{\SetFigFont{14}{16.8}{\rmdefault}{\mddefault}{\updefault}{\color[rgb]{0,0,0}$\bse'$}%
}}}}
\put(5401,-1936){\makebox(0,0)[b]{\smash{{\SetFigFont{14}{16.8}{\rmdefault}{\mddefault}{\updefault}{\color[rgb]{0,0,0}$\bsv''$}%
}}}}
\put(5401,-286){\makebox(0,0)[b]{\smash{{\SetFigFont{14}{16.8}{\rmdefault}{\mddefault}{\updefault}{\color[rgb]{0,0,0}$\bsv'$}%
}}}}
\put(6376,-961){\makebox(0,0)[b]{\smash{{\SetFigFont{14}{16.8}{\rmdefault}{\mddefault}{\updefault}{\color[rgb]{0,0,0}$\bseb$}%
}}}}
\put(601,-811){\makebox(0,0)[b]{\smash{{\SetFigFont{14}{16.8}{\rmdefault}{\mddefault}{\updefault}{\color[rgb]{0,0,0}$\bsvb$}%
}}}}
\put(1276,-811){\makebox(0,0)[b]{\smash{{\SetFigFont{14}{16.8}{\rmdefault}{\mddefault}{\updefault}{\color[rgb]{0,0,0}$\bse$}%
}}}}
\put(1726,-811){\makebox(0,0)[b]{\smash{{\SetFigFont{14}{16.8}{\rmdefault}{\mddefault}{\updefault}{\color[rgb]{0,0,0}$\bsv$}%
}}}}
\end{picture}%

%% file: e1.1.pstex_t
\begin{picture}(0,0)%
\includegraphics{e1.1.pstex}%
\end{picture}%
\setlength{\unitlength}{3947sp}%
\begingroup\makeatletter\ifx\SetFigFont\undefined%
\gdef\SetFigFont#1#2#3#4#5{%
  \reset@font\fontsize{#1}{#2pt}%
  \fontfamily{#3}\fontseries{#4}\fontshape{#5}%
  \selectfont}%
\fi\endgroup%
\begin{picture}(6171,1578)(517,-1648)
\put(3376,-286){\makebox(0,0)[b]{\smash{{\SetFigFont{17}{20.4}{\rmdefault}{\mddefault}{\updefault}{\color[rgb]{0,0,0}$\graph(\bsv)$}%
}}}}
\put(1201,-586){\makebox(0,0)[b]{\smash{{\SetFigFont{12}{14.4}{\rmdefault}{\mddefault}{\updefault}{\color[rgb]{0,0,0}$b_1$}%
}}}}
\put(1201,-1486){\makebox(0,0)[b]{\smash{{\SetFigFont{12}{14.4}{\rmdefault}{\mddefault}{\updefault}{\color[rgb]{0,0,0}$b_1$}%
}}}}
\put(2476,-1486){\makebox(0,0)[b]{\smash{{\SetFigFont{12}{14.4}{\rmdefault}{\mddefault}{\updefault}{\color[rgb]{0,0,0}$b_2$}%
}}}}
\put(2401,-586){\makebox(0,0)[b]{\smash{{\SetFigFont{12}{14.4}{\rmdefault}{\mddefault}{\updefault}{\color[rgb]{0,0,0}$b_2$}%
}}}}
\put(4576,-1036){\makebox(0,0)[b]{\smash{{\SetFigFont{12}{14.4}{\rmdefault}{\mddefault}{\updefault}{\color[rgb]{0,0,0}$b_1$}%
}}}}
\put(6376,-1036){\makebox(0,0)[b]{\smash{{\SetFigFont{12}{14.4}{\rmdefault}{\mddefault}{\updefault}{\color[rgb]{0,0,0}$b_2$}%
}}}}
\end{picture}%

%% file: 2005-41.bbl
\begin{thebibliography}

\bibitem{hb:remarks}
\textbf{H Bass}, \emph{Some remarks on group actions on trees}, Comm. Algebra 4
  (1976) 1091--1126 \MR{0419616}

\bibitem{bf:outerlimits}
\textbf{M Bestvina}, \textbf{M Feighn}, \emph{Outer limits}, preprint (1994)

\bibitem{bf:combination}
\textbf{M Bestvina}, \textbf{M Feighn}, \emph{A combination theorem for
  negatively curved groups}, J. Differential Geom. 35 (1992) 85--101
  \MR{1152226}

\bibitem{bf:bounding}
\textbf{M Bestvina}, \textbf{M Feighn}, \emph{Bounding the complexity of
  simplicial group actions on trees}, Invent. Math. 103 (1991) 449--469
  \MR{1091614}

\bibitem{bh:nonpositive}
\textbf{M\,R Bridson}, \textbf{A Haefliger}, \emph{Metric spaces of
  non-positive curvature}, Grundlehren series 319,
  Springer-Verlag, Berlin (1999) \MR{1744486}

\bibitem{bw:squares}
\textbf{M\,R Bridson}, \textbf{D\,T Wise}, \emph{{$\mathcal{VH}$} complexes,
  towers and subgroups of {$F\times F$}}, Math. Proc. Cambridge Philos. Soc.
  126 (1999) 481--497 \MR{1684244}

\bibitem{pb:splittings}
\textbf{P Brinkmann}, \emph{Splittings of mapping tori of free group
  automorphisms}, Geom. Dedicata 93 (2002) 191--203 \MR{1934698}

\bibitem{dc:book}
\textbf{D\,E Cohen}, \emph{Combinatorial group theory: a topological approach},
  London Mathematical Society Student Texts 14, Cambridge University Press,
  Cambridge (1989) \MR{1020297}

\bibitem{gad:thesis}
\textbf{G-A Diao}, \emph{Is a graph of finitely generated free groups free?
  {A}n algorithm}, PhD thesis, Rutgers University, Newark (2003)

\bibitem{fh:coherence}
\textbf{M Feighn}, \textbf{M Handel}, \emph{Mapping tori of free group
  automorphisms are coherent}, Ann. of Math. (2) 149 (1999) 1061--1077
  \MR{1709311}

\bibitem{fp:jsj}
\textbf{K Fujiwara}, \textbf{P Papasoglu}, \emph{{JSJ}-decompositions of
  finitely presented groups and complexes of groups}, to appear in GAFA

\bibitem{gmsw:hopfian}
\textbf{R Geoghegan}, \textbf{M\,L Mihalik}, \textbf{M Sapir}, \textbf{D\,T
  Wise}, \emph{Ascending {HNN} extensions of finitely generated free groups are
  {H}opfian}, Bull. London Math. Soc. 33 (2001) 292--298 \MR{1817768}

\bibitem{vg:connectedness}
\textbf{V Gerasimov}, \emph{Detecting connectedness of the boundary of a
  hyperbolic group}, preprint (1999)

\bibitem{sg:whitehead}
\textbf{S\,M Gersten}, \emph{On {W}hitehead's algorithm}, Bull. Amer. Math.
  Soc. (N.S.) 10 (1984) 281--284 \MR{733696}

\bibitem{ig:grushko}
\textbf{IA Grushko}, \emph{On generators of a free product of groups}, Matem.
  Sbornik N. S. 8 (1940) 169--182

\bibitem{ah:topology}
\textbf{A Hatcher}, \emph{Algebraic topology}, Cambridge University Press,
  Cambridge (2002) \MR{1867354}

\bibitem{jlr:algorithms}
\textbf{W Jaco}, \textbf{D Letscher}, \textbf{J\,H Rubinstein},
  \emph{Algorithms for essential surfaces in 3-manifolds}, from: ``Topology and
  geometry: commemorating SISTAG'', Contemp. Math. 314, Amer. Math. Soc.,
  Providence, RI (2002)  107--124 \MR{1941626}

\bibitem{sk:gersten}
\textbf{S Kalajd{\v{z}}ievski}, \emph{Automorphism group of a free group:
  centralizers and stabilizers}, J. Algebra 150 (1992) 435--502 \MR{1176906}

\bibitem{km:limit1}
\textbf{O Kharlampovich}, \textbf{A Myasnikov}, \emph{Irreducible affine
  varieties over a free group. {I}. {I}rreducibility of quadratic equations and
  {N}ullstellensatz}, J. Algebra 200 (1998) 472--516 \MR{1610660}

\bibitem{km:algorithms}
\textbf{O Kharlampovich}, \textbf{A Myasnikov}, \emph{Effective {JSJ}
  decompositions}, from: ``Groups, languages, algorithms'', (Borovik, editor),
  Contemp. Math. 378, Amer. Math. Soc., Providence, RI (2005)  87--212

\bibitem{ls:book}
\textbf{R\,C Lyndon}, \textbf{P\,E Schupp}, \emph{Combinatorial group theory},
  Classics in Mathematics, Springer-Verlag, Berlin (2001) \MR{1812024}

\bibitem{cfm:decision}
\textbf{C\,F Miller, III}, \emph{On group-theoretic decision problems and their
  classification}, Princeton University Press, Princeton, N.J. (1971)
  \MR{0310044}

\bibitem{rs:jsj}
\textbf{E Rips}, \textbf{Z Sela}, \emph{Cyclic splittings of finitely presented
  groups and the canonical {JSJ} decomposition}, Ann. of Math. (2) 146 (1997)
  53--109 \MR{1469317}

\bibitem{zs:tarski1}
\textbf{Z Sela}, \emph{Diophantine geometry over groups. {I}.
  {M}akanin-{R}azborov diagrams}, Publ. Math. Inst. Hautes \'Etudes Sci.
  (2001) 31--105 \MR{1863735}

\bibitem{se:trees}
\textbf{J-P Serre}, \emph{Trees}, Springer Monographs in Mathematics,
  Springer-Verlag, Berlin (2003) \MR{1954121}

\bibitem{sh:amalg}
\textbf{A Shenitzer}, \emph{Decomposition of a group with a single defining
  relation into a free product}, Proc. Amer. Math. Soc. 6 (1955) 273--279
  \MR{0069174}

\bibitem{st:folding}
\textbf{J\,R Stallings}, \emph{Topology of finite graphs}, Invent. Math. 71
  (1983) 551--565 \MR{695906}

\bibitem{js:grushko}
\textbf{J\,R Stallings}, \emph{Foldings of {$G$}-trees}, from: ``Arboreal group
  theory (Berkeley, CA, 1988)'', Math. Sci. Res. Inst. Publ. 19, Springer, New
  York (1991)  355--368 \MR{1105341}

\bibitem{swarup:hnn}
\textbf{G\,A Swarup}, \emph{Decompositions of free groups}, J. Pure Appl.
  Algebra 40 (1986) 99--102 \MR{825183}

\bibitem{jhcw:certain}
\textbf{J\,H\,C Whitehead}, \emph{On certain sets of elements in a free group},
  Proc.~London~Math.~Soc. 41 (1936) 48--56

\bibitem{jhcw:equivalent}
\textbf{J\,H\,C Whitehead}, \emph{On equivalent sets of elements in a free
  group}, Ann. of Math. (2) 37 (1936) 782--800 \MR{MR1503309}

\end{thebibliography}
